\documentclass[12pt]{amsart}
\usepackage{amssymb}
\def\opn#1#2{\def#1{\operatorname{#2}}} 
\opn\chara{char}
\opn\length{\ell}
\opn\projdim{proj\,dim}
\opn\injdim{inj\,dim}
\opn\rank{rank}
\opn\depth{depth}
\opn\grade{grade}
\opn\height{height}
\opn\embdim{emb\,dim}
\opn\codepth{codepth}
\opn\codim{codim}

\opn\Tr{Tr}
\opn\bigrank{big\,rank}
\opn\superheight{superheight}\opn\lcm{lcm}
\opn\trdeg{tr\,deg}%
\opn\reg{reg}
\opn\ini{in}
\opn\div{div}
\opn\Div{Div}
\opn\cl{cl}
\opn\Cl{Cl}
\opn\Spec{Spec}
\opn\Supp{Supp}
\opn\supp{supp}
\opn\Sing{Sing}
\opn\Ass{Ass}
\opn\Ann{Ann}
\opn\Rad{Rad}
\opn\Soc{Soc}
\opn\Ker{Ker}
\opn\Coker{Coker}
\opn\Im{Im}
\opn\Hom{Hom}
\opn\Tor{Tor}
\opn\Ext{Ext}
\opn\End{End}
\opn\Aut{Aut}
\opn\id{id}

\opn\nat{nat}
\opn\pff{pf}
\opn\Pf{Pf}
\opn\GL{GL}
\opn\SL{SL}
\opn\mod{mod}
\opn\ord{ord}
\opn\aff{aff}
\opn\con{conv}
\opn\relint{relint}
\opn\st{st}
\opn\lk{lk}
\opn\cn{cn}
\opn\core{core}
\opn\vol{vol}
\opn\gr{gr}

\def\pot#1#2{#1[\kern-0.28ex[#2]\kern-0.28ex]}

\opn\dirlim{\underrightarrow{\lim}}
\opn\invlim{\underleftarrow{\lim}}
\def\Implies{\ifmmode\Longrightarrow \else
     \unskip${}\Longrightarrow{}$\ignorespaces\fi}
\def\implies{\ifmmode\Rightarrow \else
     \unskip${}\Rightarrow{}$\ignorespaces\fi}
\def\iff{\ifmmode\Longleftrightarrow \else
     \unskip${}\Longleftrightarrow{}$\ignorespaces\fi}

\let\:=\colon
\newtheorem{Theorem}{Theorem}[section]
\newtheorem{Lemma}[Theorem]{Lemma}

\let\epsilon=\varepsilon
\let\phi=\varphi
\let\kappa=\varkappa
\def\qed{\ifhmode\textqed\fi
   \ifmmode\ifinner\quad\qedsymbol\else\dispqed\fi\fi}
\def\textqed{\unskip\nobreak\penalty50
    \hskip2em\hbox{}\nobreak\hfil\qedsymbol
    \parfillskip=0pt \finalhyphendemerits=0}
\def\dispqed{\rlap{\qquad\qedsymbol}}

\begin{document}
\title[Primitive Polynomials With Three Coefficients Given]{Existence of primitive\\ polynomials with three\\ coefficients prescribed}
\author{Donald Mills}
\address{Department of Mathematics, Southern Illinois University-Carbondale, Carbondale, IL 62901-4408}
\email{dmills@math.siu.edu}
\thanks {The author, at the time he began this project, was a Davies Fellow for the National Research Council.  He wishes to thank the NRC, and specifically the U.S. Army Research Laboratory and the U.S. Military Academy, for the use of their facilities.}

\date{\today}
\maketitle
%
%
%
%
%
%
\section{Introduction}

Let $\mathbf{F}_{q}$ denote the finite field of $q$ elements, $q=p^{r}$ for prime $p$ and positive integer $r$. A monic polynomial $f(x)=x^{n}+\sum_{i=1}^{n}f_{i}x^{n-i} \in \mathbf{F}_{q}[x]$ is called a {\it primitive polynomial\/} if it is irreducible over $\mathbf{F}_{q}$ and any of the roots of $f$ can be used to generate the multiplicative group $\mathbf{F}_{q^{n}}^{\ast}$ of $\mathbf{F}_{q^{n}}$. Equivalently, $f$ is primitive if the smallest positive integer $w$ such that $f(x) \mid x^{w}-1$ is $w=q^{n}-1$. Primitive polynomials and their roots are of interest in various applications of finite fields to coding theory and cryptography, and so it is of interest to know whether for a given $q$ and $n$ there exists a primitive polynomial of degree $n$ over $\mathbf{F}_{q}$ which may satisfy certain additional conditions. One such condition is whether there exists a primitive polynomial of degree $n$ over $\mathbf{F}_{q}$ with first coefficient $f_{1}$ prescribed, where we note that $f_{1}=-Tr(\alpha)$, $\alpha$ a root of $f$ and $Tr$ the trace function from $\mathbf{F}_{q^{n}}$ to $\mathbf{F}_{q}$. This question has been answered (see \cite{coh}, \cite{jun}), with answer as given in Theorem 1.1.

%
%

\begin{Theorem}
Let $n>1$ be an integer, and let $a \in \mathbf{F}_{q}$ be given. Then there always exists a primitive polynomial $f(x)=x^{n}+\sum_{i=1}^{n}f_{i}x^{n-i} \in \mathbf{F}_{q}[x]$ such that $a=f_{1}$ provided $(a,n) \neq (0,3)$ for $q=4$ and $(a,n) \neq (0,2)$ for arbitrary $q$.
\end{Theorem}

Cohen, Han and Mills considered the case in which there exists a primitive polynomial with $f_{1}$ and $f_{2}$ prescribed. Han \cite{han} gave the following; this result was also addressed in \cite{comi}. 

%
%

\begin{Theorem}
Let $n \geq 7$ be an integer, and let $a$, $b \in \mathbf{F}_{q}$ be given, $q$ an odd prime power. Then there always exists a primitive polynomial $f(x)=x^{n}+\sum_{i=1}^{n}f_{i}x^{n-i} \in \mathbf{F}_{q}[x]$ such that $f_{1}=a$ and $f_{2}=b$.
\end{Theorem}

Equivalently, $N_{q,n}(a,b)>0$ for all odd prime powers $q$ and all integers $n \geq 7$, where $N_{q,n}(a,b)$ is the number of primitive polynomials in $\mathbf{F}_{q}$ of degree $n$ with root $\alpha$ such that $Tr(\alpha)=a$ and $Tr(\alpha^{2})=b$, $Tr$ the trace function from $\mathbf{F}_{q^{n}}$ to $\mathbf{F}_{q}$. The case where $q=2^i$ for some $i$ is more difficult; a discussion of this case is provided in \cite{sun}.

From Theorem 1.2, we infer that the remaining cases of interest are $n=4$, $5$, and $6$. Using sieving techniques due to Cohen, Cohen and Mills \cite{comi} proved the following, with $q$ an odd prime power.

%
%

\begin{Theorem}
For all pairs $a$, $b \in \mathbf{F}_{q}$, $q$ odd, $N_{q,n}(a,b)>0$ for $n=5,6$. 
\end{Theorem}

In this paper, we generalize the above work by producing a formula in Section 2, over finite fields of suitably large characteristic, for the $k$th coefficient of an irreducible polynomial. We then use this formula to address the question of the existence of primitive polynomials with three coefficients prescribed over finite fields of characteristic at least five. The main result of the paper is given as Theorem 7.1, which states that for all finite fields of characteristic at least five, and for all $n \geq 9$, for every triplet $(f_{1},f_{2},f_{3}) \in \mathbf{F}_{q}^{3}$ there exists a primitive polynomial of degree $n$ with $x^{n-i}$ coefficient equal to $f_{i}$ for $i=1$, $2$, $3$. Progress is also made on the cases $n=7$ and $n=8$; Section $8$ is devoted to a consideration of these cases.

As the formula in Section $2$ applies to irreducible polynomials in general, and not merely primitive polynomials, the author is confident that the formula may prove important in several applications, and not only with regards to the question of existence of certain primitive polynomials.

\section{A Recursive Formula for the $k$th Coefficient of a Polynomial}

Let $f(x)=x^{n}+\sum_{i=1}^{n}(-1)^{i}f_{i}x^{n-i} \in \mathbb{F}_{q}[x]$ be given. For positive integers $k$ and $n$, $k<n$, set 

\begin{eqnarray}
W_{k,n}(x) & = & x^{q+q^2+\cdots+q^{k-2}}\sum_{i_{1}=k-1}^{n-1}x^{q^{i_{1}}}\nonumber\\ 
& + & x^{q+q^2+\cdots+q^{k-3}+q^{k-1}}\sum_{i_{2}=k}^{n-1}x^{q^{i_{2}}}\nonumber\\ 
& + & \cdots+x^{q^{n-k+1}+q^{n-k+2}+\cdots+q^{n-1}}.
\end{eqnarray}

Observe that the number of terms in $W_{k,n}$, denoted by $Z_{k,n}$, is $Z_{k,n}=\binom{n-1}{k-1}$, with $W_{1,n}(x):=1$, $W_{2,n}(x)=x^{q}+x^{q^{2}}+ \cdots +x^{q^{n-1}}$, and so forth. We have the following.

\begin{Lemma}\label{fk-form}
Let $f(x)=x^n+\sum_{i=1}^{n}(-1)^{i}f_{i}x^{n-i} \in \mathbb{F}_{q}[x]$, $p=$char$(\mathbb{F}_{q})$, denote an irreducible of degree $n$ over $\mathbb{F}_{q}$ with root $\alpha$, and let $k<n$ be any positive integer with $p \dagger k$. Then $f_{k}=\frac{1}{k}Tr(\alpha W_{k,n}(\alpha))$.
\end{Lemma}

We first prove the following technical lemma.

\begin{Lemma}\label{tech}
For positive coprime integers $k$ and $n$, $k<n$, $\frac{1}{n}\binom{n}{k}$ is integral.
\end{Lemma}

\begin{proof}
Observe that, by Legendre's identity [see for example page 67 of T. Apostol's {\it Introduction to Analytic Number Theory\/} text] we have

\begin{eqnarray}
(n-1)! & = & \prod_{p \leq n-1}p^{\sum_{m=1}^{\infty}\lfloor\frac{n-1}{p^{m}}\rfloor}
\end{eqnarray}

\noindent
for $p$ prime. Similar equations can be given for $(n-k)!$ and $k!$; note that the exponent for each $p$ has only a finite number of terms as $\lfloor\frac{n-1}{p^{m}}\rfloor = 0$ for $p^{m}>n-1$. Thinking of $\binom{n}{k}$ as 

\begin{eqnarray} 
n\frac{(n-1)(n-2)\cdots (n-k+1)}{k!},
\end{eqnarray} 

\noindent
we have by (2) that

\begin{eqnarray}
& & n\frac{(n-1)(n-2)\cdots (n-k+1)}{k!}\nonumber\\ & = & n\frac{\left(\prod_{n-k<p<n}p^{\sum_{m=1}^{\infty}\left\lfloor\frac{n-1}{p^{m}}\right\rfloor}\right)\left(\prod_{p \leq n-k}p^{\sum_{m=1}^{\infty}\left\lfloor\frac{n-1}{p^{m}}\right\rfloor-\left\lfloor\frac{n-k}{p^{m}}\right\rfloor}\right)}{\prod_{p \leq k}p^{\sum_{m=1}^{\infty}\left\lfloor\frac{k}{p^{m}}\right\rfloor}}.
\end{eqnarray}

\noindent
Thus we need to show that for each prime $p \leq k$,

\begin{eqnarray}
\sum_{m=1}^{\infty}\left\lfloor\frac{n-1}{p^{m}}\right\rfloor-\left\lfloor\frac{n-k}{p^{m}}\right\rfloor & \geq & \sum_{m=1}^{\infty}\left\lfloor\frac{k}{p^{m}}\right\rfloor.
\end{eqnarray}

If $p \mid k$ then $p \dagger (n-k)$ as $\gcd(k,n)=1$. So for each $m$, $\displaystyle \left\lfloor\frac{n-k}{p^{m}}\right\rfloor=\left\lfloor\frac{n-k-1}{p^{m}}\right\rfloor$, thus $\displaystyle \left\lfloor\frac{n-k}{p^{m}}\right\rfloor+\left\lfloor\frac{k}{p^{m}}\right\rfloor \leq \left\lfloor\frac{n-1}{p^{m}}\right\rfloor$ for each $m$ by definition of the floor function. Note that if $p \mid n$ we can make the same argument, thus the only remaining case is the one in which prime $p$ divides neither $k$ nor $n$. Since $p \dagger k$, though, we have $\displaystyle \left\lfloor\frac{k}{p^{m}}\right\rfloor=\left\lfloor\frac{k-1}{p^{m}}\right\rfloor$ and we can conclude, as before, that $\displaystyle \left\lfloor\frac{n-k}{p^{m}}\right\rfloor+\left\lfloor\frac{k}{p^{m}}\right\rfloor \leq \left\lfloor\frac{n-1}{p^{m}}\right\rfloor$ for each $m$ by definition of the floor function. Appealing to equation (5) completes the proof.
\end{proof}

\begin{proof}
We now prove Lemma \ref{fk-form}, first for the case where $\gcd(k,n)=1$. Let $\mathbb{Z}$ denote the set of integers, and let $\mathbb{Z}_{n}$ denote the set of integers modulo $n$. Observe that for each vector of the form $(a_{1},a_{2},...,a_{k}) \in \mathbb{Z}_{n}^{k}$ corresponding to the positive integer $q^{a_{1}}+q^{a_{2}}+ \cdots +q^{a_{k}}<q^{n}-1$ (assuming without loss of generality that $a_{1}<a_{2}<\cdots < a_{k}$ modulo $n$, so that the number of vectors to consider is $\binom{n}{k}$ -- let $M$ denote the set of all such vectors), it follows that as we raise $q^{a_{1}}+q^{a_{2}}+ \cdots +q^{a_{k}}$ by powers of $q$, doing our work modulo $q^{n}-1$, that the number of distinct integers formed modulo $q^{n}-1$ is exactly $n$, for $\gcd(k,n)=1$. From this observation, coupled with Lemma \ref{tech}, we conclude that $M$ can be partitioned into classes, with each class having exactly $k$ elements with $a_{1}=0$. Such an element can serve as the representative of the class.

Thus, $f_{k}$ can be written as

\begin{eqnarray}\label{eq6}
f_{k} & = & \sum_{0 \leq i_{1}<i_{2}< \cdots <i_{k}<n}\alpha^{q^{i_{1}}+q^{i_{2}}+\cdots+q^{i_{k}}} \nonumber\\
& = & \frac{1}{k}\sum_{0 \leq j_{1},j_{2},...,j_{k}<n}\alpha^{q^{j_{1}}+q^{j_{2}}+\cdots+q^{j_{k}}}
\end{eqnarray}

\noindent
where the second sum in (\ref{eq6}) amounts to $k$ copies of the first sum, the second sum having $k\binom{n}{k}=n\binom{n-1}{k-1}$ terms. Referring back to (1), observe that if one takes the trace of $\alpha W_{k,n}(\alpha)$, one obtains an expression in $\alpha$ having $n\binom{n-1}{k-1}$ terms, with each of the $ \binom{n}{k}$ members of $M$ appearing as an exponent of $\alpha$ exactly $k$ times. Thus we deduce that $f_{k}=\frac{1}{k}Tr(\alpha W_{k,n}(\alpha))$, and the first statement is proved.


For the case $1<\gcd(k,n) \leq k$ we proceed in a similar manner. Specifically, we note that the number of distinct integers modulo $q^{n}-1$ that one forms (as one raises by powers of $q$) will always be a multiple of $n/\gcd(k,n)$, say $d=\frac{ns}{\gcd(k,n)}$ for some $s$. Observe that $d$ divides $n$ as well. In applying the second expression for (\ref{eq6}), we note that each exponent in the class is found $k$ times, for a total of $dk$ terms. On the other hand, letting $L$ denote the number of exponents $q^{a_{1}}+q^{a_{2}}+ \cdots +q^{a_{k}}$ in the class with $a_{1}=0$ (without loss of generality having $a_{1} < a_{2} < \cdots < a_{k}$), applying the trace to $\alpha W_{k,n}(\alpha)$ shows that each exponent in the class appears $\frac{Ln}{d}$ times, for a total of $Ln$ terms. If $d=n$, so that $L=k$, then $dk=Ln$ and we proceed as in the first part of the proof. If $d$ is a proper divisor of $n$, then, by separating the set $\{0,1,...,n-1\}$ into equally-sized blocks of size $d$, namely into the sets $\{md,md+1,...,md+d-1\}$ for $m=0,1,...,(n/d)-1$, and by considering the coset representative $(a_{1},a_{2},...,a_{k})=(0,a_{2},...,a_{k})$ of the class, we deduce that each block of size $d$ must not only have the same number of $a_{i}'s$, but if $a_{j}=r$ for some $r$ between $0$ and $d-1$ then $a_{j}+md$ must belong to the coset representative for $m$ from $1$ to $(n/d)-1$. That is, each such class possesses a symmetry in accordance with the value of $d$. (For example, consider the case $k=6$, $n=14$ with $(a_{1},...,a_{6})=(0,1,2,7,8,9)$.) From this we conclude immediately that $L=(d/n)k$ or $dk=Ln$, as was the case for $\gcd(k,n)=1$. Arguing as above, we conclude that $f_{k}=\frac{1}{k}Tr(\alpha W_{k,n}(\alpha))$, and the lemma is proven.

\end{proof}

We now use Lemma \ref{fk-form} to prove the following.

\begin{Theorem}\label{fk-recur}
Under the conditions set forth in Lemma \ref{fk-form}, we can write $f_{k}$ as

\begin{eqnarray}
f_{k} & = & \frac{1}{k}\left(f_{k-1}Tr(\alpha)-f_{k-2}Tr(\alpha^2)+ \cdots +(-1)^{k-1}Tr(\alpha^k)\right).
\end{eqnarray}
\end{Theorem}

\begin{proof}
With $f_{0}=1$, the result is trivial for $k=1$. Observe that for any positive $k$, $$\displaystyle f_{k}=W_{k+1,n}(\alpha)+\alpha W_{k,n}(\alpha).$$ Since $$f_{k}=\frac{1}{k}Tr(\alpha W_{k,n}(\alpha))$$ by Lemma \ref{fk-form}, for all $k>1$ with $p$ not dividing $k$ we have

\begin{eqnarray*}
f_{k} & = & \frac{1}{k}Tr\left[\alpha\left(f_{k-1}-\alpha W_{k-1,n}(\alpha)\right)\right]\\
& = & \frac{1}{k}\left[f_{k-1}Tr(\alpha)-Tr\left(\alpha^{2}[f_{k-2}-\alpha W_{k-2,n}(\alpha)]\right)\right]\\
& = & \cdots\\ 
& = & \frac{1}{k}\left(f_{k-1}Tr(\alpha)-f_{k-2}Tr(\alpha^2)+ \cdots +(-1)^{k-1}Tr(\alpha^{k})\right).
\end{eqnarray*}
\end{proof}

Clearly, in order to consider the case in which $f_{1}$, $f_{2}$, and $f_{3}$ are prescribed, we will need $p \geq 5$. The character sum analysis will show that we need to restrict $k$ such that $k \leq \lfloor \frac{n}{2} \rfloor$.

\section{Character Sum Analysis}

%

To guarantee the existence of primitive polynomials over a finite field $\mathbf{F}_{q}$ with $f_{1}$, $f_{2}$, and $f_{3}$ prescribed, we will employ a character sum analysis followed by, for certain $n$, a combinatorial sieve due to Cohen (see \cite{chco}).

We first give a definition. An element $x \in \mathbf{F}_{q}$ is said to be {\it $e-$free\/} (it has also been referred to as ``no kind of $e$th power"; see for example \cite{coh2}) if, for any $y \in \mathbf{F}_{q}$ with $y^{d}=x$ for $d \mid e$, we must have $d=1$. Thus the primitive elements of $\mathbf{F}_{q}$ are those which are $(q-1)$-free, while (trivially) all elements of the field are $1$-free.

Now let $e$ denote a divisor of $q^{n}-1$, where $q$, $n$, and $a$, $b$, $c \in \mathbf{F}_{q}$ are given, and let $N(e)$ denote the number of elements $x \in \mathbf{F}_{q^{n}}$ that are $e-$free, with $Tr(x)=a$, $Tr(x^{2})=b$, and $Tr(x^{3})=c$. Further let $\omega(z)$ denote the number of prime divisors of $z$. We have the following basic lemmas.

\begin{Lemma}\label{vino}

For $\xi \in \mathbf{F}_{q^{n}}^{\ast}$, we have

\begin{eqnarray}\label{vinoequ}
\frac{\phi(e)}{e}\sum_{d \mid e}\frac{\mu(d)}{\phi(d)}\sum_{\chi_{d}}\chi^{(d)}(\xi),
\end{eqnarray}

\noindent
which equals 1 if $\xi$ is not any kind of $e$th power, and equals zero otherwise. Here $\phi$ and $\mu$ are the Euler-phi and M\"obius functions, respectively, and the inner sum runs over all $d$th-order multiplicative characters of $\mathbf{F}_{q^{n}}$. 
\end{Lemma}

\begin{Lemma}\label{basic}
For $\xi \in \mathbf{F}_{q}$ and $\psi_{t}$ an additive character of $\mathbf{F}_{q}$ for $t \in \mathbf{F}_{q}$, we have

\begin{eqnarray}
\sum_{t \in \mathbf{F}_{q}}\psi_{t}(\xi)=q
\end{eqnarray}

\noindent
if $\xi=0$. The sum equals zero otherwise.
\end{Lemma}

Using these lemmas, we may write $N(e)$ as

\begin{eqnarray}\label{exprNe}
q^{3}N(e) & = & \theta(e)\sum_{d \mid e}\frac{\mu(d)}{\phi(d)}\sum_{\chi^{(d)}}\sum_{e_{1},e_{2},e_{3} \in \mathbf{F}_{q}}S_{d,e_{1},e_{2},e_{3}}
\end{eqnarray}

\noindent
where 

\begin{eqnarray}\label{Sequ}
S_{d,e_{1},e_{2},e_{3}} & = & \sum_{\xi \in \mathbf{F}_{q^{n}}^{\ast}}\psi\left(Tr(e_{1}\xi+e_{2}\xi^{2}+e_{3}\xi^{3})-e_{1}a-e_{2}b-e_{3}c\right)\chi^{(d)}(\xi),
\end{eqnarray}

\noindent
$\chi^{(d)}$ runs over all $d$th-order multiplicative characters of $\mathbf{F}_{q^{n}}$, $\psi$ is the canonical $\mathbf{F}_{q}$-additive character, $\theta(e)=\phi(e)/e$, and $Tr$ is the trace map from $\mathbf{F}_{q^{n}}$ to $\mathbf{F}_{q}$.

Observe that when $e=q^{n}-1$, $N:=N(q^{n}-1)$ is the value whose positivity we wish to determine; note as well that the value of $N(e)$ depends only on the distinct prime factors of $e$. With these observations in hand, we say that divisors $e_{1}$, ..., $e_{r}$, $r \geq 1$, of $e$ are {\it complementary divisors of $e$ with common divisor $d$\/} if the set of distinct prime divisors of lcm$\{e_{1},...,e_{r}\}$ is the same as that of $e$, and, for any pair $(i,j)$ with $1 \leq i \neq j \leq r$, the set of distinct prime divisors of $\gcd(e_{i},e_{j})$ is that of $d$. When $r=1$, we have $e_{1}=d=e$.

With these notions in hand, we arrive at the following sieve inequality, proved in \cite{chco}.

\begin{Theorem}
Let $q$ be a prime power and $n \geq 1$ an integer. Let $e_{1}$, ..., $e_{r}$, $r \geq 1$ be complementary divisors of $e \mid q^{n}-1$ with common divisor $d$. Then, with $N(e)$ defined as above, we have

\begin{eqnarray}\label{generalsieve}
N(e) & \geq & \left[\sum_{i=1}^{r}N(e_{i})\right]-(r-1)N(d).
\end{eqnarray}
\end{Theorem}

Thus it suffices to guarantee

\begin{eqnarray}\label{applysieve}
\left[\sum_{i=1}^{r}N(e_{i})\right]-(r-1)N(d) & > & 0.
\end{eqnarray}

Before using the sieve, we must obtain bounds for $N$, depending upon the values of $a$, $b$, and $c$. First, we note that the following lemma will prove useful \cite{han}.

\begin{Lemma}\label{twisted}
Let $\chi$ denote a $d$th order multiplicative character and $\psi$ an additive character of $\mathbf{F}_{q}$. Let $f(x)$, $g(x) \in \mathbf{F}_{q}[x]$ be polynomials of degree $m$, $r$ respectively. If $\gcd(m,d)=\gcd(r,q)=1$, then

\begin{eqnarray*}
\left|\sum_{c \in \mathbf{F}_{q}}\chi(f(c))\psi(g(c))\right| & \leq & (m+r-1)\sqrt{q}.
\end{eqnarray*} 
\end{Lemma}

Of course, $S_{1,0,0,0}=q^{n}-1$. We have the following.

\begin{Theorem}\label{Nbound}
We have

\begin{eqnarray}\label{BoundforN}
q^{3}N & \geq & \theta(q^{n}-1)\{q^{n}-1+T_{1}-\sum_{i=2}^{8}|T_{i}|\},
\end{eqnarray}

\noindent
where the $T_{i}$, $i=1$, ..., $8$ are defined below.
\end{Theorem}

\begin{proof}
Our work is separated into the following cases, based upon the values of $d$ and the $e_{j}$:

\begin{enumerate}
\item{$d=1$, $e_{j}=0$ for all $j$ (addressed above).}
\item{$d=1$, $e_{j} \neq 0$ for exactly one $i$.}
\item{$d=1$, $e_{j} \neq 0$ for exactly two $i$.}
\item{$d=1$, $e_{1}e_{2}e_{3} \neq 0$.}
\item{$d>1$, distinguishing as to whether $d \mid Q:=\frac{q^{n}-1}{q-1}$.}
\end{enumerate}

For Case $2$, we have the following subcases.\\

(2a) $e_{1} \neq 0$. The sum to consider is 

\begin{eqnarray*}
T_{1} & = & \sum_{e_{1} \in \mathbf{F}_{q}^{\ast}}S_{1,e_{1},0,0}\nonumber\\
& = & \sum_{e_{1} \in \mathbf{F}_{q}^{\ast}}\sum_{\xi \in \mathbf{F}_{q^{n}}^{\ast}}\psi(Tr(e_{1}\xi)-e_{1}a)\nonumber\\
& = & \sum_{e_{1} \in \mathbf{F}_{q}^{\ast}}\psi(-e_{1}a)\left[\sum_{\xi \in \mathbf{F}_{q^{n}}}\psi(Tr(e_{1}\xi))-1\right]\nonumber\\
& = & \sum_{e_{1} \in \mathbf{F}_{q}^{\ast}}\psi(-e_{1}a)\left[\sum_{\xi \in \mathbf{F}_{q^{n}}}\psi(Tr(\xi))-1\right].\nonumber\\
\end{eqnarray*}

Thus,

\begin{center}
\[ T_{1}= \left\{ \begin{array}{ll}
                   1-q & \mbox{if $a=0$}\\
                   1   & \mbox{if $a \neq 0$.}
                 \end{array} 
         \right. \]
\end{center}

(2b) $e_{2} \neq 0$. The sum to consider is

\begin{eqnarray*}
T_{2} & = & \sum_{e_{2} \in \mathbf{F}_{q}^{\ast}}S_{1,0,e_{2},0}\nonumber\\
& = & \sum_{e_{2} \in \mathbf{F}_{q}^{\ast}}\sum_{\xi \in \mathbf{F}_{q^{n}}^{\ast}}\psi(Tr(e_{2}\xi^{2})-e_{2}b).\nonumber\\
\end{eqnarray*}

From \cite{han}, we have

\begin{center}
\[ |T_{2}| \leq \left\{ \begin{array}{ll}
                        (q-1)(\sqrt{q^{n}}+1)          & \mbox{if $b=0$}\\
                        (\sqrt{q}+1)(\sqrt{q^{n}}+1)   & \mbox{if $b \neq 0$.}
                        \end{array} 
                \right. \]
\end{center}

(2c) $e_{3} \neq 0$. We consider

\begin{eqnarray*}
T_{3} & = & \sum_{e_{3} \in \mathbf{F}_{q}^{\ast}}S_{1,0,0,e_{3}}\nonumber\\
& = & \sum_{e_{3} \in \mathbf{F}_{q}^{\ast}}\sum_{\xi \in \mathbf{F}_{q^{n}}^{\ast}}\psi(Tr(e_{3}\xi^{3})-e_{3}c)\nonumber\\
& = & \sum_{e_{3} \in \mathbf{F}_{q}^{\ast}}\psi(-e_{3}c)\sum_{\xi \in \mathbf{F}_{q^{n}}^{\ast}}\lambda(e_{3}\xi^{3}),
\end{eqnarray*}

\noindent
where $\lambda(X)=\psi(Tr(X))$ for all $X \in \mathbf{F}_{q^{n}}$. We divide the work here into two subcases.\\

(2c1) $q \equiv 2 \pmod 3$. Thus $\gcd(3,q-1)=1$, and so we may write $T_{3}$ as

\begin{eqnarray*}
T_{3} & = & \sum_{e_{3} \in \mathbf{F}_{q}^{\ast}}\psi(-e_{3}c)\sum_{\xi \in \mathbf{F}_{q^{n}}^{\ast}}\lambda(e_{3}\xi^{3})\nonumber\\
& = & \sum_{e_{3} \in \mathbf{F}_{q}^{\ast}}\psi(-e_{3}^{3}c)\sum_{\xi \in \mathbf{F}_{q^{n}}^{\ast}}\lambda((e_{3}\xi)^{3})\nonumber\\
& = & \sum_{e_{3} \in \mathbf{F}_{q}^{\ast}}\psi(-e_{3}^{3}c)\sum_{\xi \in \mathbf{F}_{q^{n}}^{\ast}}\lambda(\xi^{3})
\end{eqnarray*}

Thus, by the Weil bound \cite{ln} we have

\begin{center}
\[ |T_{3}| \leq \left\{ \begin{array}{ll}
                        (q-1)(2\sqrt{q^{n}}+1)           & \mbox{if $c=0$}\\
                        (2\sqrt{q}+1)(2\sqrt{q^{n}}+1)   & \mbox{if $c \neq 0$.}
                        \end{array} 
                \right. \]
\end{center}

(2c2) $q \equiv 1 \pmod 3$. Let $\alpha$ denote a fixed cubic nonresidue in $\mathbf{F}_{q}^{\ast}$, and let $C$ denote the set of cubic residues in $\mathbf{F}_{q}$. Observe that $C \cup C\alpha \cup C\alpha^{2}=\mathbf{F}_{q}^{\ast}$. We have

\begin{eqnarray*}
T_{3} & = & \frac{1}{3}\left(\sum_{i=0}^{2}\sum_{e_{3} \in \mathbf{F}_{q}^{\ast}}\sum_{\xi \in \mathbf{F}_{q^{n}}^{\ast}}\psi(Tr(e_{3}^{3}\alpha^{i}\xi^{3})-e_{3}^{3}\alpha^{i}c)\right)\nonumber\\
& = & \frac{1}{3}\left(\sum_{i=0}^{2}\sum_{e_{3} \in \mathbf{F}_{q}^{\ast}}\psi(-e_{3}^{3}\alpha^{i}c)\sum_{\xi \in \mathbf{F}_{q^{n}}^{\ast}}\lambda(\alpha^{i}\xi^{3})\right).\nonumber\\
\end{eqnarray*}

Thus, we again have

\begin{center}
\[ |T_{3}| \leq \left\{ \begin{array}{ll}
                        (q-1)(2\sqrt{q^{n}}+1)           & \mbox{if $c=0$}\\
                        (2\sqrt{q}+1)(2\sqrt{q^{n}}+1)   & \mbox{if $c \neq 0$.}
                        \end{array} 
                \right. \]
\end{center}

For Case 3, we also have three subcases to address.\\

(3a) $e_{1}e_{2} \neq 0$. The sum in question is

\begin{eqnarray*}
T_{4} & = & \sum_{e_{1},e_{2} \in \mathbf{F}_{q}^{\ast}}\sum_{\xi \in \mathbf{F}_{q^{n}}^{\ast}}\psi(Tr(e_{1}\xi+e_{2}\xi^{2})-e_{1}a-e_{2}b),\nonumber
\end{eqnarray*}

\noindent
whose modulus is bounded from above by (see \cite{han})

\begin{center}
\[ |T_{4}| \leq \left\{ \begin{array}{lll}
                        (q-1)^{2}(\sqrt{q^{n}}+1)           & \mbox{if $a=b=0$}\\
                        (q-1)(\sqrt{q^{n}}+1)               & \mbox{if $a \neq 0$, $b=0$}\\
                        (q-1)(\sqrt{q}+1)(\sqrt{q^{n}}+1)   & \mbox{if $b \neq 0$.}
                        \end{array} 
                \right. \]
\end{center}

(3b) $e_{1}e_{3} \neq 0$. We consider

\begin{eqnarray*}
T_{5} & = & \sum_{e_{1},e_{3} \in \mathbf{F}_{q}^{\ast}}\sum_{\xi \in \mathbf{F}_{q^{n}}^{\ast}}\psi(Tr(e_{1}\xi+e_{3}\xi^{3})-e_{1}a-e_{3}c)\nonumber\\
& = & \sum_{e,e_{1} \in \mathbf{F}_{q}^{\ast}}\sum_{\xi \in \mathbf{F}_{q^{n}}^{\ast}}\psi(Tr(e_{1}\xi+e(e_{1}\xi)^{3})-e_{1}a-ee_{1}^{3}c)\nonumber\\
& = & \sum_{e \in \mathbf{F}_{q}^{\ast}}\sum_{\xi \in \mathbf{F}_{q^{n}}^{\ast}}\lambda(\xi+e\xi^{3})\sum_{e_{1} \in \mathbf{F}_{q}^{\ast}}\psi(-e_{1}a-ee_{1}^{3}c)\nonumber\\
\end{eqnarray*}

\noindent
where $ee_{1}^{3}=e_{3}$ and $\lambda$ has the same meaning as above. Thus

\begin{center}
\[ |T_{5}| \leq \left\{ \begin{array}{lll}
                        (q-1)^{2}(2\sqrt{q^{n}}+1)           & \mbox{if $a=c=0$}\\
                        (q-1)(2\sqrt{q^{n}}+1)               & \mbox{if $a \neq 0$, $c=0$}\\
                        (q-1)(2\sqrt{q}+1)(2\sqrt{q^{n}}+1)  & \mbox{if $c \neq 0$.}
                        \end{array} 
                \right. \]
\end{center}

(3c) $e_{2}e_{3} \neq 0$. Let $\alpha$ denote a fixed quadratic nonresidue in $\mathbf{F}_{q}$. Consider the sum

\begin{eqnarray*}
T_{6} & = & \sum_{e_{2},e_{3} \in \mathbf{F}_{q}^{\ast}}\sum_{\xi \in \mathbf{F}_{q^{n}}^{\ast}}\psi(Tr(e_{2}\xi^{2}+e_{3}\xi^{3})-e_{2}b-e_{3}c)\nonumber\\
& = & \frac{1}{2}\left(\sum_{i=0}^{1}\sum_{e_{2},e_{3} \in \mathbf{F}_{q}^{\ast}}\sum_{\xi \in \mathbf{F}_{q^{n}}^{\ast}}\psi(Tr(e_{2}^{2}\alpha^{i}\xi^{2}+e_{3}\xi^{3})-e_{2}^{2}\alpha^{i}b-e_{3}c)\right)\nonumber\\
& = & \frac{1}{2}\left(\sum_{i=0}^{1}\sum_{e,e_{2} \in \mathbf{F}_{q}^{\ast}}\sum_{\xi \in \mathbf{F}_{q^{n}}^{\ast}}\psi(Tr(\alpha^{i}(e_{2}\xi)^{2}+e(e_{2}\xi)^{3})-\alpha^{i}e_{2}^{2}b-ee_{2}^{3}c)\right)\nonumber\\
& = & \frac{1}{2}\left(\sum_{i=0}^{1}\sum_{e \in \mathbf{F}_{q}^{\ast}}\sum_{\xi \in \mathbf{F}_{q^{n}}^{\ast}}\lambda(\alpha^{i}\xi^{2}+e\xi^{3})\sum_{e_{2} \in \mathbf{F}_{q}^{\ast}}\psi(-\alpha^{i}e_{2}^{2}b-ee_{2}^{3}c)\right),
\end{eqnarray*}

\noindent
where $ee_{2}^{3}=e_{3}$. Thus

\begin{center}
\[ |T_{6}| \leq \left\{ \begin{array}{lll}
                        (q-1)^{2}(2\sqrt{q^{n}}+1)                       & \mbox{if $b=c=0$}\\
                        (q-1)(\sqrt{q}+1)(2\sqrt{q^{n}}+1)               & \mbox{if $b \neq 0$, $c=0$}\\
                        (q-1)(2\sqrt{q}+1)(2\sqrt{q^{n}}+1)              & \mbox{if $c \neq 0$.}
                        \end{array} 
                \right. \]
\end{center}

(4) $e_{1}e_{2}e_{3} \neq 0$. We have

\begin{eqnarray*}
T_{7} & = & \sum_{e_{1},e_{2},e_{3} \in \mathbf{F}_{q}^{\ast}}\sum_{\xi \in \mathbf{F}_{q^{n}}^{\ast}}\psi(Tr(e_{1}\xi+e_{2}\xi^{2}+e_{3}\xi^{3})-e_{1}a-e_{2}b-e_{3}c).\nonumber\\
\end{eqnarray*}

\noindent
Set $ee_{1}^{2}=e_{2}$ and $ge_{1}^{3}=e_{3}$ now and proceed as before to obtain

\begin{eqnarray*}
T_{7} & = & \sum_{e,g \in \mathbf{F}_{q}^{\ast}}\sum_{\xi \in \mathbf{F}_{q^{n}}^{\ast}}\lambda(\xi+e\xi^{2}+g\xi^{3})\sum_{e_{1} \in \mathbf{F}_{q}^{\ast}}\psi(-e_{1}a-ee_{1}^{2}b-ge_{1}^{3}c).\nonumber\\
\end{eqnarray*}

Thus,

\begin{center}
\[ |T_{7}| \leq \left\{ \begin{array}{llll}
                        (q-1)^{3}(2\sqrt{q^{n}}+1)                           & \mbox{if $a=b=c=0$}\\
                        (q-1)^{2}(2\sqrt{q^{n}}+1)                           & \mbox{if $a \neq 0$, $b=c=0$}\\
                        (q-1)^{2}(\sqrt{q}+1)(2\sqrt{q^{n}}+1)               & \mbox{if $b \neq 0$, $c=0$}\\
                        (q-1)^{2}(2\sqrt{q}+1)(2\sqrt{q^{n}}+1)              & \mbox{if $c \neq 0$.}
                        \end{array} 
                \right. \]
\end{center}

(5) $d>1$. Here we use the fact that the $\mathbf{F}_{q^{n}}$-multiplicative character $\chi^{(d)}$, applied to $\mathbf{F}_{q}$, is trivial iff $d \mid Q$. The sum to consider is

\begin{eqnarray}\label{T8}
T_{8} & = & \sum_{1<d \mid q^{n}-1}\frac{\mu(d)}{\phi(d)}\sum_{\chi^{(d)}}\sum_{e_{1},e_{2},e_{3} \in \mathbf{F}_{q}}S_{d,e_{1},e_{2},e_{3}}.
\end{eqnarray}

\noindent
We will consider (\ref{T8}) according to the values of the $e_{i}$'s, specifically as to whether a certain $e_{i}=0$. We shall also separate our results according to whether $d \mid Q$. Note that $S_{d,0,0,0}=0$. Proceeding in the same manner as above, we have (with $\alpha$ a quadratic nonresidue of $\mathbf{F}_{q}$ in (\ref{equS2}), and $ee_{1}^{2}=e_{2}$ in (\ref{equS3}))

\begin{eqnarray}\label{equS1}
\sum_{e_{1} \in \mathbf{F}_{q}^{\ast}}S_{d,e_{1},0,0} & = & \sum_{e_{1} \in \mathbf{F}_{q}^{\ast}}\chi^{(d)}(e_{1}^{-1})\psi(-(e_{1}-1)a)S_{d,1,0,0},
\end{eqnarray}

\begin{eqnarray}\label{equS2}
\sum_{e_{2} \in \mathbf{F}_{q}^{\ast}}S_{d,0,e_{2},0} & = & \frac{1}{2}\left(\sum_{i=0}^{1}\sum_{e_{2} \in \mathbf{F}_{q}^{\ast}}\chi^{(d)}(e_{2}^{-1})\psi(-(e_{2}^{2}-1)\alpha^{i}b)S_{d,0,\alpha^{i},0}\right),
\end{eqnarray}

\noindent
and

\begin{eqnarray}\label{equS3}
\sum_{e_{1},e_{2} \in \mathbf{F}_{q}^{\ast}}S_{d,e_{1},e_{2},0} & = & \sum_{e,e_{1} \in \mathbf{F}_{q}^{\ast}}\chi^{(d)}(e_{1}^{-1})\psi(-(e_{1}-1)a-(e_{1}^{2}-1)eb)S_{d,1,e,0}.
\end{eqnarray}

\noindent
Further, we have

\begin{eqnarray}\label{equS4p1}
\sum_{e_{3} \in \mathbf{F}_{q}^{\ast}}S_{d,0,0,e_{3}} & = & \sum_{e_{3} \in \mathbf{F}_{q}^{\ast}}\chi^{(d)}(e_{3}^{-1})\psi(-(e_{3}^{3}-1)c)S_{d,0,0,1}
\end{eqnarray}

\noindent
for $q \equiv 2 \pmod 3$, while, with $q \equiv 1 \pmod 3$ and $\alpha$ a fixed cubic nonresidue in $\mathbf{F}_{q}$, we have

\begin{eqnarray}\label{equS4p2}
\sum_{e_{3} \in \mathbf{F}_{q}^{\ast}}S_{d,0,0,e_{3}} & = & \frac{1}{3}\left(\sum_{i=0}^{2}\sum_{e_{3} \in \mathbf{F}_{q}^{\ast}}\chi^{(d)}(e_{3}^{-1})\psi(-(e_{3}^{3}-1)\alpha^{i}c)S_{d,0,0,\alpha^{i}}\right).
\end{eqnarray}

\noindent
With $ee_{1}^{3}=e_{3}$, we have

\begin{eqnarray}\label{equS5}
\sum_{e_{1},e_{3} \in \mathbf{F}_{q}^{\ast}}S_{d,e_{1},0,e_{3}} & = & \sum_{e,e_{1} \in \mathbf{F}_{q}^{\ast}}\chi^{(d)}(e_{1}^{-1})\psi(-(e_{1}-1)a-(e_{1}^{3}-1)ec)S_{d,1,0,e},
\end{eqnarray}

\noindent
while, with $\alpha$ a fixed quadratic nonresidue in $\mathbf{F}_{q}$ and $ee_{2}^{3}=e_{3}$, we have

\begin{eqnarray}\label{equS6}
\sum_{e,e_{2} \in \mathbf{F}_{q}^{\ast}}S_{d,0,e_{2},e_{3}} & = & \frac{1}{2}\left(\sum_{i=0}^{1}\sum_{e,e_{2} \in \mathbf{F}_{q}^{\ast}}\chi^{(d)}(e_{2}^{-1})U(e,e_{2},\alpha^{i})S_{d,0,\alpha^{i},e}\right),
\end{eqnarray}

\noindent
where $U(e,e_{2},\alpha^{i})=\psi(-(e_{2}^{2}-1)\alpha^{i}b-(e_{2}^{3}-1)ec)$. Finally, with $ee_{1}^{2}=e_{2}$ and $ge_{1}^{3}=e_{3}$, we have

\begin{eqnarray}
\sum_{e_{1},e_{2},e_{3} \in \mathbf{F}_{q}^{\ast}}S_{d,e_{1},e_{2},e_{3}} & = & \sum_{e_{1},e,g \in \mathbf{F}_{q}^{\ast}}\chi^{(d)}(e_{1}^{-1})V(e_{1},e,g)S_{d,1,e,g},
\end{eqnarray}

\noindent
where $V(e_{1},e,g)=\psi(-(e_{1}-1)a-(e_{1}^{2}-1)eb-(e_{1}^{3}-1)gc)$.

Putting it all together, we have the following bounds for $|T_{8}|$, depending upon the values of $a$, $b$, and $c$.

For $a=b=c=0$, we have

\begin{eqnarray}\label{T8bnd1}
|T_{8}| & \leq & (2^{\omega(Q)}-1)[(q-1)(3q^{2}+2q+1)]\sqrt{q^{n}}.
\end{eqnarray}

For $a \neq 0$, $b=c=0$ we have

\begin{eqnarray}\label{T8bnd2}
|T_{8}| & \leq & (2^{\omega(Q)}-1)[1+10(q-1)+6(q-1)^{2}]\sqrt{q^{n}}\nonumber\\
& + & (2^{\omega(q^{n}-1)}-2^{\omega(Q)})[1+5(q-1)+3(q-1)^{2}]\sqrt{q^{n+1}},
\end{eqnarray}

\noindent
while for $b \neq 0$ and $a=c=0$ we have

\begin{eqnarray}\label{T8bnd3}
|T_{8}| & \leq & (2^{\omega(Q)}-1)[(5q-3)(\sqrt{q}+1)+4(q-1)+3(\sqrt{q}+2)(q-1)^{2}]\sqrt{q^{n}}\nonumber\\
& + & (2^{\omega(q^{n}-1)}-2^{\omega(Q)})[4+10(q-1)+6(q-1)^{2}]\sqrt{q^{n+1}},
\end{eqnarray}

\noindent
and for $c \neq 0$ with $a=b=0$ we have

\begin{eqnarray}\label{T8bnd4}
|T_{8}| & \leq & (2^{\omega(Q)}-1)[(6q-3)(2\sqrt{q}+1)+3(q-1)+(6\sqrt{q}+5)(q-1)^{2}]\sqrt{q^{n}}\nonumber\\
& + & (2^{\omega(q^{n}-1)}-2^{\omega(Q)})[9+18(q-1)+9(q-1)^{2}]\sqrt{q^{n+1}}.
\end{eqnarray}

For $ab \neq 0$ with $c=0$, we have

\begin{eqnarray}\label{T8bnd5}
|T_{8}| & \leq & (2^{\omega(Q)}-1)[2\sqrt{q}+3+(5\sqrt{q}+11)(q-1)+(3\sqrt{q}+3)(q-1)^{2}]\sqrt{q^{n}}\nonumber\\
& + & (2^{\omega(q^{n}-1)}-2^{\omega(Q)})[5+13(q-1)+6(q-1)^{2}]\sqrt{q^{n+1}},
\end{eqnarray}

\noindent
while for $ac \neq 0$ and $b=0$ we have

\begin{eqnarray}\label{T8bnd6}
|T_{8}| & \leq & (2^{\omega(Q)}-1)[6\sqrt{q}+4+(12\sqrt{q}+10)(q-1)+(6\sqrt{q}+3)(q-1)^{2}]\sqrt{q^{n}}\nonumber\\
& + & (2^{\omega(q^{n}-1)}-2^{\omega(Q)})[10+20(q-1)+9(q-1)^{2}]\sqrt{q^{n+1}}.
\end{eqnarray}

\noindent
When $bc \neq 0$ and $a=0$, we have

\begin{eqnarray}\label{T8bnd7}
|T_{8}| & \leq & (2^{\omega(Q)}-1)[8\sqrt{q}+5+(14\sqrt{q}+9)(q-1)+(6\sqrt{q}+3)(q-1)^{2}]\sqrt{q^{n}}\nonumber\\
& + & (2^{\omega(q^{n}-1)}-2^{\omega(Q)})[13+22(q-1)+9(q-1)^{2}]\sqrt{q^{n+1}}.
\end{eqnarray}

\noindent
Finally, for $abc \neq 0$ we have

\begin{eqnarray}\label{T8bnd7}
|T_{8}| & \leq & (2^{\omega(Q)}-1)[8\sqrt{q}+6+(14\sqrt{q}+8)(q-1)+(6\sqrt{q}+3)(q-1)^{2}]\sqrt{q^{n}}\nonumber\\
& + & (2^{\omega(q^{n}-1)}-2^{\omega(Q)})[14+22(q-1)+9(q-1)^{2}]\sqrt{q^{n+1}}.
\end{eqnarray}

Putting this all together, we obtain (\ref{BoundforN}). This completes the proof.
\end{proof}

This completes the main portion of our character sum analysis. The next section is devoted to using Theorem \ref{Nbound} to make statements of the following type: ``For a given triple $(a,b,c)$, if $q^{A(n)} \geq B(q,n)$ for some functions $A$ and $B$, then $N>0$ for all fields $\mathbf{F}_{q}$ having characteristic at least $5$, and with $n \geq 7$." We will use these bounds to ensure that $N>0$ for $n \geq 13$, then move to a sieving process to resolve, as best we can, the cases $7 \leq n \leq 12$.

\section{Bounds That Ensure $N>0$}

We separate the work into two cases, depending upon whether $(a,b,c)=(0,0,0)$.

\begin{enumerate}
\item{$(a,b,c)=(0,0,0)$. Observe that (abbreviating $N_{q,n}(a,b,c)$ with $N$)

\begin{eqnarray}\label{Nzeros}
q^{3}N & \geq & \theta(q^{n}-1)[q^{n}-q-(q-1)(3\sqrt{q^{n}}+2)]\nonumber\\
& - & \theta(q^{n}-1)[(q-1)^{2}(5\sqrt{q^{n}}+3)+(q-1)^{3}(2\sqrt{q^{n}}+1)]\nonumber\\
& - & \theta(q^{n}-1)[(2^{\omega(Q)}-1)(3q^{\frac{n}{2}+3}-q^{\frac{n}{2}+2}-q^{\frac{n}{2}+1}-q^{\frac{n}{2}})].
\end{eqnarray}

\noindent
Thus, to ensure that $N>0$ it suffices to ensure that

\begin{eqnarray}\label{condNzeros}
q^{n}-q-(q-1)(3\sqrt{q^{n}}+2)-(q-1)^{2}(5\sqrt{q^{n}}+3)-(q-1)^{3}(2\sqrt{q^{n}}+1)\nonumber\\
-(2^{\omega(Q)}-1)(3q^{\frac{n}{2}+3}-q^{\frac{n}{2}+2}-q^{\frac{n}{2}+1}-q^{\frac{n}{2}}) > 0,
\end{eqnarray}

\noindent
or, by grouping the terms in (\ref{condNzeros}) with minus signs in front of them,

\begin{eqnarray}\label{cNzeros}
(3 \times 2^{\omega(Q)}+2)q^{\frac{n}{2}+3}+6q^{\frac{n}{2}+2}+10q^{\frac{n}{2}+1}\nonumber\\
+6q^{\frac{n}{2}}+q^{3}+3q^{2}+6q+3<q^{n}.
\end{eqnarray}

\noindent
Thus, we want

\begin{eqnarray}\label{penultNzeros}
3 \times 2^{\omega(Q)}+2+\frac{6}{q}+\frac{10}{q^{2}}+\frac{6}{q^{3}}+\frac{1}{q^{\frac{n}{2}}}+\frac{3}{q^{\frac{n}{2}+1}}+\frac{6}{q^{\frac{n}{2}+2}}+\frac{3}{q^{\frac{n}{2}+3}}<q^{\frac{n}{2}-3},
\end{eqnarray}

\noindent
or, since $q \geq 5$ and $n \geq 7$,

\begin{eqnarray}\label{ultNzeros}
3 \times 2^{\omega(Q)} + 3.655 < q^{\frac{n}{2}-3}.
\end{eqnarray}
}
\item{$(a,b,c) \neq (0,0,0)$. We present an analysis of each of the seven cases in which $(a,b,c) \neq (0,0,0)$; these are handled in a manner like that of the case $a=b=c=0$, and in each case we will give bounds to ensure $N>0$ for said case, as we did for (\ref{ultNzeros}). 

For $a \neq 0$, $b=c=0$ we have

\begin{eqnarray}\label{Nanotzero}
q^{3}N & \geq & \theta(q^{n}-1)[q^{n}-(6\sqrt{q^{n}}+4)(q-1)-(4\sqrt{q^{n}}+2)(q-1)^{2}]\nonumber\\
& - & \theta(q^{n}-1)(2^{\omega(Q)}-1)[1+10(q-1)+6(q-1)^{2}]\sqrt{q^{n}}\nonumber\\
& - & \theta(q^{n}-1)(2^{\omega(q^{n}-1)}-2^{\omega(Q)})[1+5(q-1)+3(q-1)^{2}]\sqrt{q^{n+1}}.
\end{eqnarray}

\noindent
Thus to ensure $N>0$, we proceed as in the all-zero case to conclude that we want

\begin{eqnarray}\label{cNaneqzero}
2^{\omega(q^{n}-1)}(3q^{\frac{n+5}{2}})+2^{\omega(Q)}(6q^{\frac{n+4}{2}}+q^{\frac{n+3}{2}}+q^{\frac{n+1}{2}})\nonumber\\
+4q^{\frac{n+4}{2}}+8q^{\frac{n+2}{2}}+7q^{\frac{n}{2}}+2q^{2}+4q+2<q^{n},
\end{eqnarray}

\noindent
or, replacing $2^{\omega(Q)}$ with $2^{\omega(q^{n}-1)}$ and recalling that $q \geq 5$, $n \geq 7$, we have, after dividing through on both sides by $q^{\frac{n+5}{2}}$ and then setting $q=5$ and $n=7$ where appropriate, rounding up to the nearest thousandth in our work,

\begin{eqnarray}\label{ultNanotzero}
(5.924)2^{\omega(q^{n}-1)}+2.635 < q^{\frac{n-5}{2}}.
\end{eqnarray}

\vspace{0.3in}

For $b \neq 0$, $a=c=0$, we have

\begin{eqnarray}\label{Nbnotzero}
q^{3}N & \geq & \theta(q^{n}-1)[q^{n}-(q+(\sqrt{q}+1)(\sqrt{q^{n}}+1))]\nonumber\\
& - & \theta(q^{n}-1)[2\sqrt{q^{n}}+1+(\sqrt{q}+1)(3\sqrt{q^{n}}+2)](q-1)\nonumber\\
& - & \theta(q^{n}-1)(2\sqrt{q^{n}}+1)(\sqrt{q}+2)(q-1)^{2}\nonumber\\
& - & \theta(q^{n}-1)(2^{\omega(Q)}-1)[(5q-3)(\sqrt{q}+1)+4(q-1)+3(\sqrt{q}+2)(q-1)^{2}]\sqrt{q^{n}}\nonumber\\
& - & \theta(q^{n}-1)(2^{\omega(q^{n}-1)}-2^{\omega(Q)})[4+10(q-1)+6(q-1)^{2}]\sqrt{q^{n+1}}.
\end{eqnarray}

\noindent
Thus, to ensure that $N>0$, we want

\begin{eqnarray}\label{cNbneqzero}
2^{\omega(q^{n}-1)}(6q^{\frac{n+5}{2}})+2^{\omega(Q)}(3q^{\frac{n+5}{2}}+6q^{\frac{n+4}{2}}+2q^{\frac{n+3}{2}})\nonumber\\
+2q^{\frac{n+5}{2}}+4q^{\frac{n+4}{2}}+4q^{\frac{n+3}{2}}+8q^{\frac{n+2}{2}}+3q^{\frac{n+1}{2}}+6q^{\frac{n}{2}}\nonumber\\
+q^{\frac{5}{2}}+2q^{2}+2q^{\frac{3}{2}}+4q+2q^{\frac{1}{2}}+3<q^{n},
\end{eqnarray}

\noindent
or, replacing $2^{\omega(Q)}$ with $2^{\omega(q^{n}-1)}$ and recalling that $q \geq 5$, $n \geq 7$, we have

\begin{eqnarray}\label{ultNbnotzero}
(12.083)2^{\omega(q^{n}-1)}+5.542 < q^{\frac{n-5}{2}}.
\end{eqnarray}

\vspace{0.3in}

For $c \neq 0$, $a=b=0$, the inequality to consider is

\begin{eqnarray}\label{Ncnotzero}
q^{3}N & \geq & \theta(q^{n}-1)[q^{n}-(q+(2\sqrt{q}+1)(2\sqrt{q^{n}}+1))]\nonumber\\
& - & \theta(q^{n}-1)[2(2\sqrt{q}+1)(2\sqrt{q^{n}}+1)+\sqrt{q^{n}}+1](q-1)\nonumber\\
& - & \theta(q^{n}-1)[\sqrt{q^{n}}+1+(2\sqrt{q}+1)(2\sqrt{q^{n}}+1)](q-1)^{2}\nonumber\\
& - & \theta(q^{n}-1)(2^{\omega(Q)}-1)[(6q-3)(2\sqrt{q}+1)+3(q-1)+(6\sqrt{q}+5)(q-1)^{2}]\sqrt{q^{n}}\nonumber\\
& - & \theta(q^{n}-1)(2^{\omega(q^{n}-1)}-2^{\omega(Q)})[9+18(q-1)+9(q-1)^{2}]\sqrt{q^{n+1}}.
\end{eqnarray}

\noindent
Thus, to ensure that $N>0$, it suffices to have (replacing $\omega(Q)$ with $\omega(q^{n}-1)$)

\begin{eqnarray}\label{cNcneqzero}
2^{\omega(q^{n}-1)}(15q^{\frac{n+5}{2}}+5q^{\frac{n+4}{2}})+4q^{\frac{n+5}{2}}+3q^{\frac{n+4}{2}}\nonumber\\
+8q^{\frac{n+3}{2}}+6q^{\frac{n+2}{2}}+8q^{\frac{n+1}{2}}+6q^{\frac{n}{2}}+2q^{\frac{5}{2}}+2q^{2}\nonumber\\
+4q^{\frac{3}{2}}+4q+4\sqrt{q}+3<q^{n},
\end{eqnarray}

\noindent
or

\begin{eqnarray}\label{ultNcnotzero}
(15+\sqrt{5})2^{\omega(q^{n}-1)}+7.921<q^{\frac{n-5}{2}}.
\end{eqnarray}

\vspace{0.3in}

Now consider the case $ab \neq 0$ with $c=0$. We have

\begin{eqnarray}\label{Nabnotzero}
q^{3}N & \geq & \theta(q^{n}-1)[q^{n}-(\sqrt{q}+1)(\sqrt{q^{n}}+1))]\nonumber\\
& - & \theta(q^{n}-1)[(\sqrt{q}+1)(3\sqrt{q^{n}}+2)+4\sqrt{q^{n}}+2](q-1)\nonumber\\
& - & \theta(q^{n}-1)(\sqrt{q}+1)(2\sqrt{q^{n}}+1)(q-1)^{2}\nonumber\\
& - & \theta(q^{n}-1)(2^{\omega(Q)}-1)[2\sqrt{q}+3+(5\sqrt{q}+11)(q-1)+(3\sqrt{q}+3)(q-1)^{2}]\sqrt{q^{n}}\nonumber\\
& - & \theta(q^{n}-1)(2^{\omega(q^{n}-1)}-2^{\omega(Q)})[5+13(q-1)+6(q-1)^{2}]\sqrt{q^{n+1}}.
\end{eqnarray}

\noindent
Thus, to ensure that $N>0$, it suffices to have (replacing $\omega(Q)$ with $\omega(q^{n}-1)$)

\begin{eqnarray}\label{cNabneqzero}
2^{\omega(q^{n}-1)}\left(9q^{\frac{n+5}{2}}+3q^{\frac{n+4}{2}}+q^{\frac{n+3}{2}}+5q^{\frac{n+2}{2}}+2q^{\frac{n+1}{2}}\right)\nonumber\\
+2q^{\frac{n+5}{2}}+2q^{\frac{n+4}{2}}+4q^{\frac{n+3}{2}}+7q^{\frac{n+2}{2}}+3q^{\frac{n+1}{2}}+8q^{\frac{n}{2}}\nonumber\\
+q^{\frac{5}{2}}+q^{2}+2q^{\frac{3}{2}}+4q+2q^{\frac{1}{2}}+2<q^{n},
\end{eqnarray}

\noindent
or

\begin{eqnarray}\label{ultNabnotzero}
(11.069)2^{\omega(q^{n}-1)}+4.592<q^{\frac{n-5}{2}}.
\end{eqnarray}

\vspace{0.3in}

Now consider the case $ac \neq 0$ with $b=0$. We have

\begin{eqnarray}\label{Nacnotzero}
q^{3}N & \geq & \theta(q^{n}-1)[q^{n}-(2\sqrt{q}+1)(2\sqrt{q^{n}}+1))]\nonumber\\
& - & \theta(q^{n}-1)[\sqrt{q^{n}}+1+(2\sqrt{q^{n}}+1)(2\sqrt{q}+1)](2q-2)\nonumber\\
& - & \theta(q^{n}-1)(2\sqrt{q}+1)(2\sqrt{q^{n}}+1)(q-1)^{2}\nonumber\\
& - & \theta(q^{n}-1)(2^{\omega(Q)}-1)[6\sqrt{q}+4+(12\sqrt{q}+10)(q-1)+(6\sqrt{q}+3)(q-1)^{2}]\sqrt{q^{n}}\nonumber\\
& - & \theta(q^{n}-1)(2^{\omega(q^{n}-1)}-2^{\omega(Q)})[10+20(q-1)+9(q-1)^{2}]\sqrt{q^{n+1}}.
\end{eqnarray}

\noindent
Thus, to ensure that $N>0$, it suffices to have

\begin{eqnarray}\label{cNacneqzero}
2^{\omega(q^{n}-1)}\left(15q^{\frac{n+5}{2}}+3q^{\frac{n+4}{2}}+2q^{\frac{n+3}{2}}+4q^{\frac{n+2}{2}}+q^{\frac{n+1}{2}}\right)\nonumber\\
+4q^{\frac{n+5}{2}}+2q^{\frac{n+4}{2}}+8q^{\frac{n+3}{2}}+6q^{\frac{n+2}{2}}+8q^{\frac{n+1}{2}}+7q^{\frac{n}{2}}\nonumber\\
+2q^{\frac{5}{2}}+q^{2}+4q^{\frac{3}{2}}+4q+4q^{\frac{1}{2}}+2<q^{n},
\end{eqnarray}

\noindent
or

\begin{eqnarray}\label{ultNacnotzero}
(17.140)2^{\omega(q^{n}-1)}+7.490<q^{\frac{n-5}{2}}.
\end{eqnarray}

\vspace{0.3in}

Now consider $bc \neq 0$ with $a=0$. We have

\begin{eqnarray}\label{Nbcnotzero}
q^{3}N & \geq & \theta(q^{n}-1)[q^{n}-(q+(\sqrt{q}+1)(\sqrt{q^{n}}+1)+(2\sqrt{q}+1)(2\sqrt{q^{n}}+1))]\nonumber\\
& - & \theta(q^{n}-1)[(\sqrt{q}+1)(\sqrt{q^{n}}+1)+2(2\sqrt{q}+1)(2\sqrt{q^{n}}+1)](q-1)\nonumber\\
& - & \theta(q^{n}-1)(2\sqrt{q}+1)(2\sqrt{q^{n}}+1)(q-1)^{2}\nonumber\\
& - & \theta(q^{n}-1)(2^{\omega(Q)}-1)[8\sqrt{q}+5+(14\sqrt{q}+9)(q-1)+(6\sqrt{q}+3)(q-1)^{2}]\sqrt{q^{n}}\nonumber\\
& - & \theta(q^{n}-1)(2^{\omega(q^{n}-1)}-2^{\omega(Q)})[13+22(q-1)+9(q-1)^{2}]\sqrt{q^{n+1}}.
\end{eqnarray}

\noindent
Thus it suffices to have

\begin{eqnarray}\label{cNbcneqzero}
2^{\omega(q^{n}-1)}\left(15q^{\frac{n+5}{2}}+3q^{\frac{n+4}{2}}+6q^{\frac{n+3}{2}}+3q^{\frac{n+2}{2}}\right)\nonumber\\
+4q^{\frac{n+5}{2}}+2q^{\frac{n+4}{2}}+9q^{\frac{n+3}{2}}+5q^{\frac{n+2}{2}}+9q^{\frac{n+1}{2}}+6q^{\frac{n}{2}}\nonumber\\
+2q^{\frac{5}{2}}+q^{2}+5q^{\frac{3}{2}}+4q+5q^{\frac{1}{2}}+3<q^{n},
\end{eqnarray}

\noindent
or

\begin{eqnarray}\label{ultNbcnotzero}
(17.810)2^{\omega(q^{n}-1)}+7.624<q^{\frac{n-5}{2}}.
\end{eqnarray}

\vspace{0.3in}

Finally, consider $abc \neq 0$. We have

\begin{eqnarray}\label{Nabcnotzero}
q^{3}N & \geq & \theta(q^{n}-1)[q^{n}-((\sqrt{q}+1)(\sqrt{q^{n}}+1)+(2\sqrt{q}+1)(2\sqrt{q^{n}}+1))]\nonumber\\
& - & \theta(q^{n}-1)[(\sqrt{q}+1)(\sqrt{q^{n}}+1)+2(2\sqrt{q}+1)(2\sqrt{q^{n}}+1)](q-1)\nonumber\\
& - & \theta(q^{n}-1)(2\sqrt{q}+1)(2\sqrt{q^{n}}+1)(q-1)^{2}\nonumber\\
& - & \theta(q^{n}-1)(2^{\omega(Q)}-1)[8\sqrt{q}+6+(14\sqrt{q}+8)(q-1)+(6\sqrt{q}+3)(q-1)^{2}]\sqrt{q^{n}}\nonumber\\
& - & \theta(q^{n}-1)(2^{\omega(q^{n}-1)}-2^{\omega(Q)})[14+22(q-1)+9(q-1)^{2}]\sqrt{q^{n+1}}.
\end{eqnarray}

\noindent
Thus it suffices to have

\begin{eqnarray}\label{cNabcneqzero}
2^{\omega(q^{n}-1)}\left(15q^{\frac{n+5}{2}}+3q^{\frac{n+4}{2}}+6q^{\frac{n+3}{2}}+2q^{\frac{n+2}{2}}+q^{\frac{n+1}{2}}+q^{\frac{n}{2}}\right)\nonumber\\
+4q^{\frac{n+5}{2}}+2q^{\frac{n+4}{2}}+9q^{\frac{n+3}{2}}+5q^{\frac{n+2}{2}}+9q^{\frac{n+1}{2}}+5q^{\frac{n}{2}}\nonumber\\
+2q^{\frac{5}{2}}+q^{2}+5q^{\frac{3}{2}}+3q+5q^{\frac{1}{2}}+3<q^{n},
\end{eqnarray}

\noindent
or

\begin{eqnarray}\label{ultNabcnotzero}
(17.779)2^{\omega(q^{n}-1)}+7.606<q^{\frac{n-5}{2}}.
\end{eqnarray}
}
\end{enumerate}

\vspace{0.3in}

Of all the inequalities given for $(a,b,c) \neq (0,0,0)$, (\ref{ultNbcnotzero}) is the most stringent, and thus we will use this inequality, along with (\ref{ultNzeros}), to resolve the existence question in the following section for $n \geq 13$.

\section{The Case $n \geq 13$}

As in the previous section, we separate our work into two cases, according to whether $(a,b,c)=(0,0,0)$.

\vspace{0.3in}

\begin{enumerate}
\item{$(a,b,c)=(0,0,0)$. Refer to (\ref{ultNzeros}). Observe that we can strengthen this inequality to read

\begin{eqnarray}\label{strongultNzeros}
2^{\omega(Q)+2.735}<q^{\frac{n-6}{2}}
\end{eqnarray}

\noindent
or

\begin{eqnarray}\label{strong2ultNzeros}
q^{n}>\left(2^{\omega(Q)+2.735}\right)^{u_{0}(n)}
\end{eqnarray}

\noindent
where $u_{0}(n)=\frac{2n}{n-6}$ for $n \geq 7$. Consider (\ref{strong2ultNzeros}) with $n \geq 13$, so that $u_{0}(n) \leq \frac{26}{7}$. If $\omega(Q) \geq 19$ then

\begin{eqnarray}\label{ngeq13_1}
Q & \geq & A_{19} \times 2^{6.149(\omega(Q)-19)}\nonumber\\
& > & 2^{\frac{26}{7}(\omega(Q)+2.735)}\nonumber\\
& \geq & 2^{u_{0}(n)(\omega(Q)+2.735)},
\end{eqnarray}

\noindent
where $A_{19}$ is the product of the first $19$ primes. So when $\omega(Q) \geq 19$, (\ref{strong2ultNzeros}) holds. If $\omega(Q) \leq 18$ and 

\begin{eqnarray}\label{ngeq13_2}
q^{n}> 2^{\frac{(26)(20.735)}{7}}
\end{eqnarray} 

\noindent
then again (\ref{strong2ultNzeros}) holds. The $q^{n}$ values to check directly, that is, those which do not satisfy (\ref{ngeq13_2}), are: $5^{n}$, $13 \leq n \leq 33$; $7^{n}$, $13 \leq n \leq 27$; $11^{n}$, $13 \leq n \leq 22$; $13^{n}$, $13 \leq n \leq 20$; $17^{n}$ and $19^{n}$, $13 \leq n \leq 18$; $23^{n}$, $13 \leq n \leq 17$; $25^{n}$, $13 \leq n \leq 16$; $29^{n}$ and $31^{n}$, $13 \leq n \leq 15$; $37^{n}$, $41^{n}$, and $43^{n}$, $13 \leq n \leq 14$; and $47^{13}$, $49^{13}$, $53^{13}$, and $59^{13}$. All of the possible exceptions listed, however, satisfy (\ref{strong2ultNzeros}), and thus $N_{q,n}(0,0,0)>0$ for char$(\mathbf{F}_{q}) \geq 5$ and $n \geq 13$. \vspace{0.3in}
}
\item{$(a,b,c) \neq (0,0,0)$. Observe that we can strengthen (\ref{ultNbcnotzero}) to read

\begin{eqnarray}\label{strongultNbcnotzero}
2^{\omega(q^{n}-1)+4.669}<q^{\frac{n-5}{2}}
\end{eqnarray}

\noindent
or

\begin{eqnarray}\label{strong2ultNbcnotzero}
q^{n}>\left(2^{\omega(q^{n}-1)+4.669}\right)^{u_{1}(n)}
\end{eqnarray}

\noindent
where $u_{1}(n)=\frac{2n}{n-5}$ for $n \geq 7$. Consider (\ref{strong2ultNbcnotzero}) with $n \geq 13$, so that $u_{1}(n) \leq 3.25$. We have

\begin{eqnarray}\label{ngeq13_3}
q^{n} & > & Q(q-1)\nonumber\\
& > & Q(6 \times 2^{u_{1}(n)(\omega(q-1)-4)})\nonumber\\
& > & 2^{2.58+u_{1}(n)(\omega(q-1)-4)}Q.
\end{eqnarray}

\noindent
Thus if 

\begin{eqnarray}\label{ngeq13_4}
Q & > & \left(2^{\omega(Q)+4.669}\right)^{u_{1}(n)}2^{4u_{1}(n)-2.58},
\end{eqnarray}

\noindent
or, more stringently, as $u_{1}(n) \leq 3.25$,

\begin{eqnarray}\label{ngeq13_5}
Q > 2^{25.595+3.25\omega(Q)},
\end{eqnarray}

\noindent
then (\ref{strong2ultNbcnotzero}) is satisfied. Recall that for $n \geq 13$, if $\omega(Q) \geq 19$ then $Q>2^{\frac{26}{7}(\omega(Q)+2.735)}$. Thus we ask when 

\begin{eqnarray}\label{ngeq13_6}
\frac{26}{7}(\omega(Q)+2.735) & > & 25.595+3.25\omega(Q);
\end{eqnarray}

\noindent
the answer to that is that (\ref{ngeq13_6}) is satisfied for $\omega(Q) \geq 34$. Thus, $N>0$ for char$(\mathbf{F}_{q}) \geq 5$ and $n \geq 13$ with $\omega(Q) \geq 34$.

For $\omega(Q) \leq 33$, we want

\begin{eqnarray}\label{ngeq13_7}
q^{n-1} & > & \frac{2^{4u_{1}(n)}}{6}\left(2^{\omega(Q)+4.669}\right)^{u_{1}(n)},
\end{eqnarray}

\noindent
for then, as $q-1>2^{u_{1}(n)(\omega(q-1)-4)}$, it would follow that (\ref{strong2ultNbcnotzero}) is satisfied. Setting $u_{1}(n)=3.25$ and $\omega(Q)=33$ in (\ref{ngeq13_7}), it follows that we want to satisfy

\begin{eqnarray}\label{ngeq13_8}
q^{n-1} & > & 2^{132.840}.
\end{eqnarray}

We write a computer program to check whether all pairs $(q,n)$ ($q$ not a power of $2$ or $3$) that do not satisfy (\ref{ngeq13_8}), satisfy (\ref{strong2ultNbcnotzero}) nonetheless. We find that all such pairs $(q,n)$ do satisfy (\ref{strong2ultNbcnotzero}), and thus $N>0$ for char$(\mathbf{F}_{q}) \geq 5$ and $n \geq 13$.
}
\end{enumerate}

\section{Sieve Inequalities for the Three-Coefficient Problem}

We will use (\ref{applysieve}), in conjunction with the bounds given for $N$, to resolve the primitive polynomial existence question for $9 \leq n \leq 12$, and to come close to a resolution of said problem for $n=7$, $8$.

We first consider the case $a=b=c=0$. Note here that, based upon our work in bounding $N_{q,n}(0,0,0)$, and in reference to (\ref{applysieve}), we only need to work with divisors of $Q$. In particular, note that for a divisor $m$ of $Q$ we have

\begin{eqnarray}\label{genboundallzeros}
q^{3}N(m) & \geq & \theta(m)\{q^{n}-P(q,n)-(2^{\omega(m)}-1)R(q,n)\}
\end{eqnarray}

\noindent
where $\theta(m)=\phi(m)/m$,

\begin{eqnarray}\label{Pequallzeros}
P(q,n) & = & q+(q-1)(3\sqrt{q^{n}}+2)\nonumber\\
& + &(q-1)^{2}(5\sqrt{q^{n}}+3)+(q-1)^{3}(2\sqrt{q^{n}}+1),
\end{eqnarray}

\noindent
and

\begin{eqnarray}\label{Requallzeros}
R(q,n) & = & (q-1)(3q^{2}+2q+1)\sqrt{q^{n}}.
\end{eqnarray}

Observe first that

\begin{eqnarray}\label{Rboundallzeros}
R(q,n) & = & 3q^{\frac{n+6}{2}}-q^{\frac{n+4}{2}}-q^{\frac{n+2}{2}}-q^{\frac{n}{2}}\nonumber\\
& < & 3q^{\frac{n+6}{2}}.
\end{eqnarray}

\noindent
Further, after some arithmetic we find that

\begin{eqnarray}\label{Pboundallzeros}
P(q,n) & = & 2q^{\frac{n+6}{2}}-q^{\frac{n+4}{2}}-q^{\frac{n+2}{2}}+q^{3}\nonumber\\
& < & 2q^{\frac{n+6}{2}}
\end{eqnarray}

\noindent
for all prime powers $q$ with $n \geq 7$. Thus,

\begin{eqnarray}\label{revgenboundallzeros}
q^{3}N(m) & > & \theta(m)\{q^{n}-q^{\frac{n+6}{2}}(3 \times 2^{\omega(m)}-1)\}.
\end{eqnarray}

\noindent
In particular, for a set of complementary divisors $e_{1}$, ..., $e_{r}$ with common divisor $d$, we have

\begin{eqnarray}\label{rev2genboundallzeros}
\frac{q^{3}N(d)\theta}{\theta(d)} & > & \theta\{q^{n}-q^{\frac{n+6}{2}}(3 \times 2^{\omega(d)}-1)\}
\end{eqnarray}

\noindent
where $\theta:=-(r-1)\theta(d)+\sum_{i=1}^{r}\theta(e_{i})$. Here we need $\theta>0$. Now write (\ref{applysieve}) as

\begin{eqnarray}\label{apply2sieve}
\sum_{i=1}^{r}[N(e_{i})-\frac{\theta(e_{i})}{\theta(d)}N(d)]+\frac{\theta}{\theta(d)}N(d) > 0
\end{eqnarray}

\noindent
and apply (\ref{rev2genboundallzeros}), as well as

\begin{eqnarray}\label{absoluteboundallzeros}
q^{3}\left|N(e_{i})-\frac{\theta(e_{i})}{\theta(d)}N(d)\right| & \leq & 3q^{\frac{n+6}{2}}\theta(e_{i})(2^{\omega(e_{i})}-2^{\omega(d)})
\end{eqnarray}

\noindent
for each $i$, where (\ref{absoluteboundallzeros}) follows from the estimates of the character sums given earlier, as applied to those divisors of $e_{i}$ that are not involved in $N(d)$. Thus, using (\ref{rev2genboundallzeros}) and (\ref{absoluteboundallzeros}), we want

\begin{eqnarray}\label{finalsieve1allzeros}
q^{\frac{n-6}{2}} & > & \frac{3\sum_{i=1}^{r}\theta(e_{i})(2^{\omega(e_{i})}-2^{\omega(d)})}{\theta}+3 \times 2^{\omega(d)}-1
\end{eqnarray}

\noindent
in order to ensure that $N>0$. If one chooses complementary divisors such that $2^{\omega(e_{i})}-2^{\omega(d)}=2^{\omega(d)}$ for each $i$, (\ref{finalsieve1allzeros}) becomes

\begin{eqnarray}\label{finalsieve2allzeros}
q^{\frac{n-6}{2}} & > & \frac{3 \times 2^{\omega(d)}(2\theta + (r-1)\theta(d))}{\theta}-1.
\end{eqnarray}

We will use (\ref{finalsieve1allzeros}) and (\ref{finalsieve2allzeros}) in the next two sections for the cases $7 \leq n \leq 12$ and $a=b=c=0$.

We obtain the sieve inequalities for $(a,b,c) \neq (0,0,0)$ in much the same way we obtained (\ref{finalsieve1allzeros}) and (\ref{finalsieve2allzeros}). First, we consider the case $a \neq 0$, $b=c=0$. For $m$ a divisor of $q^{n}-1$, we have

\begin{eqnarray}\label{genboundaneqzero}
q^{3}N(m) \geq \theta(m)\{q^{n}-P(q,n)-(2^{\omega(m)}-1)R(q,n)\}
\end{eqnarray}

\noindent
where 

\begin{eqnarray}\label{Pboundaneqzero}
P(q,n) & = & (6\sqrt{q^{n}}+4)(q-1)+(4\sqrt{q^{n}}+2)(q-1)^{2}\nonumber\\
& < & 4q^{\frac{n+4}{2}},
\end{eqnarray}

\begin{eqnarray}\label{Rboundaneqzero}
R(q,n) & = & \left[1+10(q-1)+6(q-1)^{2}+(1+5(q-1)+3(q-1)^{2})\sqrt{q}\right]\sqrt{q^{n}}\nonumber\\
& \leq & \frac{21}{4}q^{\frac{n+5}{2}}
\end{eqnarray}

\noindent
for $q \geq 5$, and we use $2^{\omega(m)}-1$ in place of $2^{\omega(\gcd(m,Q))}-1$ and $2^{\omega(m)}-2^{\omega(\gcd(m,Q))}$. Arguing as we did for the all-zeros case, it is a straightforward matter to conclude that we want

\begin{eqnarray}\label{finalsieve1aneqzero}
q^{\frac{n-5}{2}} & > & \frac{21}{4}\left[\frac{\sum_{i=1}^{r}\theta(e_{i})(2^{\omega(e_{i})}-2^{\omega(d)})}{\theta}+2^{\omega(d)}-1\right]+\frac{4}{\sqrt{q}}
\end{eqnarray}

\noindent
or, for a choice of complementary divisors such that $2^{\omega(e_{i})}-2^{\omega(d)}=2^{\omega(d)}$ for each $i$,

\begin{eqnarray}\label{finalsieve2aneqzero}
q^{\frac{n-5}{2}} & > & \frac{(5.25)2^{\omega(d)}(2\theta+(r-1)\theta(d))}{\theta}+\frac{4}{\sqrt{q}}-\frac{21}{4}.
\end{eqnarray}

As the inequalities for the other cases are obtained in like manner, we list only the final results below, with the proofs left to the reader. For each of these cases, only the general sieve inequality is given, as the sieve inequality produced for the situation in which $2^{\omega(e_{i})}-2^{\omega(d)}=2^{\omega(d)}$ for each $i$ is easily obtained from the general expression.

For $b \neq 0$, $a=c=0$ we want

\begin{eqnarray}\label{finalsieve1bneqzero}
q^{\frac{n-5}{2}} & > & \left(9+\frac{6}{\sqrt{q}}\right)\left[\frac{\sum_{i=1}^{r}\theta(e_{i})(2^{\omega(e_{i})}-2^{\omega(d)})}{\theta}+2^{\omega(d)}-1\right]+2+\frac{4}{\sqrt{q}}.
\end{eqnarray}

For $c \neq 0$, $a=b=0$ we want

\begin{eqnarray}\label{finalsieve1cneqzero}
q^{\frac{n-5}{2}} & > & \left(15+\frac{5}{\sqrt{q}}\right)\left[\frac{\sum_{i=1}^{r}\theta(e_{i})(2^{\omega(e_{i})}-2^{\omega(d)})}{\theta}+2^{\omega(d)}-1\right]+4+\frac{3}{\sqrt{q}}.
\end{eqnarray}

For $ab \neq 0$, $c=0$ we want

\begin{eqnarray}\label{finalsieve1abneqzero}
q^{\frac{n-5}{2}} & > & \left(9+\frac{4}{\sqrt{q}}\right)\left[\frac{\sum_{i=1}^{r}\theta(e_{i})(2^{\omega(e_{i})}-2^{\omega(d)})}{\theta}+2^{\omega(d)}-1\right]+2+\frac{3}{\sqrt{q}}.
\end{eqnarray}

For $ac \neq 0$, $b=0$ we want

\begin{eqnarray}\label{finalsieve1acneqzero}
q^{\frac{n-5}{2}} & > & \left(15+\frac{5}{\sqrt{q}}\right)\left[\frac{\sum_{i=1}^{r}\theta(e_{i})(2^{\omega(e_{i})}-2^{\omega(d)})}{\theta}+2^{\omega(d)}-1\right]+4+\frac{3}{\sqrt{q}}.
\end{eqnarray}

For $bc \neq 0$, $a=0$ (and also for $abc \neq 0$) we want

\begin{eqnarray}\label{finalsieve1bcneqzero}
q^{\frac{n-5}{2}} & > & \left(15+\frac{19}{3\sqrt{q}}\right)\left[\frac{\sum_{i=1}^{r}\theta(e_{i})(2^{\omega(e_{i})}-2^{\omega(d)})}{\theta}+2^{\omega(d)}-1\right]+4+\frac{3}{\sqrt{q}}.
\end{eqnarray}

Of these inequalities, (\ref{finalsieve1bcneqzero}) is the most restrictive, and thus we shall use this inequality in the sections to follow.

\section{The Cases $9 \leq n \leq 12$}

We shall proceed in this section in descending order, beginning with $n=12$. For each section, we shall proceed in like manner to the method given in \cite{comi}. Specifically, we shall begin by using (\ref{ultNzeros}) or (\ref{ultNbcnotzero}), as appropriate, to say that $N>0$ for all $\omega(Q) \geq \omega_{0}$ or $\omega(q^{n}-1) \geq \omega_{1}$, again as appropriate, where $\omega_{0}$ and $\omega_{1}$ are determined by (\ref{ultNzeros}) or (\ref{ultNbcnotzero}), respectively. Then, we shall use (\ref{finalsieve1allzeros}) or (\ref{finalsieve1bcneqzero}), again as appropriate, to improve the results obtained using (\ref{ultNzeros}) or (\ref{ultNbcnotzero}). In this stage, we shall, for each value of $\omega$, determine a ``worst-case scenario" value $t_{\omega}$, that is, we shall be able to say for the given $\omega$ that, if the prime power $q$ in question is such that $\omega(Q)=\omega$ or $\omega(q^{n}-1)=\omega$, as appropriate, and $q>t_{\omega}$ then $N>0$ for said value $q$ and given $n$. In each case, we build a table which allows us to make such conclusions. In the first column of each table, the value of $\omega$ is given. In the second column, we give the minimum value $q_{0}$ such that $\omega$ can equal the prescribed value. In the third column, we present $t_{\omega}$, determined via the appropriate sieve inequality. (The value $t_{\omega}$ is the right-hand side, or RHS, of the sieve inequality, raised to the appropriate power in order to compare directly with $q_{0}$.) After these two steps, we will have a set of prime power values that have not yet been eliminated. These ``possible exceptions" will then be eliminated, either via the appropriate sieve inequality, or by means of direct verification, that is, we will use a computer to eliminate the prime power in question.

\begin{enumerate}
\item{$n=12$. First, consider the all-zero case. From (\ref{ultNzeros}), we determine that $N>0$ for $\omega(Q) \geq 16$. Use of (\ref{finalsieve1allzeros}) improves this to $N>0$ for $\omega(Q) \geq 10$, as Table \ref{sievetable_12_zeros} shows. In Table \ref{sievetable_12_zeros}, we use $e_{1}=d=2$ for $\omega=1$, while using $e_{1}=d=2$ and $e_{2}=6$ for $\omega=2$. For $\omega \geq 3$, we use complementary divisors such that $2^{\omega(e_{i})}-2^{\omega(d)}=2^{\omega(d)}$ for each $i$. For $\omega=3$, $d=2$, while for $\omega \geq 4$, we use $d=6$. There are no possible exceptions, as $\omega((5^{12}-1)/4)=6$. Thus, $N_{q,12}(0,0,0)>0$.


\begin{table}[htbp]\label{sievetable_12_zeros}
\centering
\begin{tabular}{||c|c|c||}\hline
$\omega(Q)$ & $q_{0}$ & cube root of RHS \\ \hline
1 & 0.50 & 1.71 \\ \hline
2 & 0.86 & 2.22 \\ \hline
3 & 1.16 & 2.88 \\ \hline
4 & 1.47 & 3.46 \\ \hline
5 & 1.89 & 4.03 \\ \hline
6 & 2.43 & 4.59 \\ \hline
7 & 3.19 & 5.12 \\ \hline
8 & 4.21 & 5.66 \\ \hline
9 & 5.64 & 6.20 \\ \hline
10 & 7.70 & 6.72 \\ \hline
\end{tabular}
\caption{Sieving Table for Case $n=12$, $a=b=c=0$}
\end{table}

For the case $(a,b,c) \neq (0,0,0)$, we deduce via (\ref{ultNbcnotzero}) that $N>0$ for $\omega(q^{12}-1) \geq 18$. Use of (\ref{finalsieve1bcneqzero}) improves this to $N>0$ for $\omega(q^{12}-1) \geq 12$, as indicated in Table 2. The construction of this table, with regards to complementary divisors, is the same as that for Table \ref{sievetable_12_zeros}. The possible exceptions here are $q=5$ and $q=7$. While the latter succumbs to the sieve for $d=2$, with $2^{\omega(e_{i})}-2^{\omega(d)}=2^{\omega(d)}$ for each $i$, the former must be checked directly, and we do so, with an affirmative outcome. Thus $N_{q,12}(a,b,c)>0$.


\begin{table}[htbp]\label{s_table_12_notallzero}
\centering
\begin{tabular}{||c|c|c||}\hline
$\omega(q^{12}-1)$ & $q_{0}$ & $(2/7)$-root of RHS \\ \hline
1 & 1.10 & 2.59 \\ \hline
2 & 1.18 & 3.36 \\ \hline
3 & 1.33 & 4.23 \\ \hline
4 & 1.56 & 4.95 \\ \hline
5 & 1.91 & 5.62 \\ \hline
6 & 2.36 & 6.25 \\ \hline
7 & 2.99 & 6.82 \\ \hline
8 & 3.82 & 7.39 \\ \hline
9 & 4.96 & 7.94 \\ \hline
10 & 6.57 & 8.46 \\ \hline
11 & 8.75 & 8.99 \\ \hline
12 & 11.82 & 9.50 \\ \hline
\end{tabular}
\caption{Sieving Table for Case $n=12$, $(a,b,c) \neq (0,0,0)$}
\end{table}}
\item{$n=11$. Again, we look at the all-zero case first. Observe first that prime $p$ divides $Q$ if and only if $p=11$ or $p \equiv 1 \pmod{22}$ (refer to page $26$ of \cite{koblitz}, for example). Using this, we determine from (\ref{ultNzeros}) that $N>0$ for $\omega(Q) \geq 5$. Use of (\ref{finalsieve1allzeros}) reflects this, as Table 3 shows. In Table 3, we use $e_{1}=d=11$ for $\omega=1$, while using $e_{1}=d=11$ and $e_{2}=253=(11)(23)$ for $\omega=2$. For $\omega \geq 3$, we use complementary divisors such that $2^{\omega(e_{i})}-2^{\omega(d)}=2^{\omega(d)}$ for each $i$. For $\omega=3$, $d=11$, while for $\omega \geq 4$, we use $d=253$. There are no possible exceptions. Thus, $N_{q,11}(0,0,0)>0$.


\begin{table}[htbp]\label{sievetable_11_zeros}
\centering
\begin{tabular}{||c|c|c||}\hline
$\omega(Q)$ & $q_{0}$ & $(2/5)$-root of RHS \\ \hline
1 & 1.00 & 1.90 \\ \hline
2 & 1.57 & 2.61 \\ \hline
3 & 2.52 & 3.13 \\ \hline
4 & 4.03 & 4.16 \\ \hline
5 & 6.93 & 4.70 \\ \hline
\end{tabular}
\caption{Sieving Table for Case $n=11$, $a=b=c=0$}
\end{table}

For the case $(a,b,c) \neq (0,0,0)$, we use (\ref{ultNbcnotzero}) to conclude that $N>0$ for $\omega(q^{11}-1) \geq 20$. Use of (\ref{finalsieve1bcneqzero}) improves this to $N>0$ for $\omega(q^{11}-1) \geq 12$, as indicated in Table 4. The construction of this table, with regards to complementary divisors, is the same as that for Table 2. The possible exception here is $q=7$, but this prime power satisfies (\ref{ultNbcnotzero}). Thus $N_{q,11}(a,b,c)>0$.


\begin{table}[htbp]\label{s_table_11_notallzero}
\centering
\begin{tabular}{||c|c|c||}\hline
$\omega(q^{11}-1)$ & $q_{0}$ & cube root of RHS \\ \hline
1 & 1.11 & 3.03 \\ \hline
2 & 1.19 & 4.10 \\ \hline
3 & 1.37 & 5.38 \\ \hline
4 & 1.63 & 6.44 \\ \hline
5 & 2.02 & 7.47 \\ \hline
6 & 2.55 & 8.45 \\ \hline
7 & 3.30 & 9.37 \\ \hline
8 & 4.32 & 10.28 \\ \hline
9 & 5.74 & 11.17 \\ \hline
10 & 7.80 & 12.02 \\ \hline
11 & 10.65 & 12.90 \\ \hline
12 & 14.79 & 13.77 \\ \hline
\end{tabular}
\caption{Sieving Table for Case $n=11$, $(a,b,c) \neq (0,0,0)$}
\end{table}}
\item{$n=10$. From (\ref{ultNzeros}), we determine that $N_{q,10}(0,0,0)>0$ for $\omega(Q) \geq 25$. Use of (\ref{finalsieve1allzeros}) improves this to $N>0$ for $\omega(Q) \geq 12$, as Table 5 shows. The table is built in the same manner as Table 1. The only possible exception is $q=5$, which, when we check directly, we find that $N_{5,10}(0,0,0)>0$. Thus, $N_{q,10}(0,0,0)>0$.


\begin{table}[htbp]\label{sievetable_10_zeros}
\centering
\begin{tabular}{||c|c|c||}\hline
$\omega(Q)$ & $q_{0}$ & square root of RHS \\ \hline
1 & 0.50 & 2.24 \\ \hline
2 & 0.88 & 3.32 \\ \hline
3 & 1.23 & 4.88 \\ \hline
4 & 1.63 & 6.42 \\ \hline
5 & 2.21 & 8.09 \\ \hline
6 & 3.01 & 9.83 \\ \hline
7 & 4.18 & 11.60 \\ \hline
8 & 5.85 & 13.48 \\ \hline
9 & 8.35 & 15.44 \\ \hline
10 & 12.19 & 17.41 \\ \hline
11 & 17.91 & 19.54 \\ \hline
12 & 26.80 & 21.72 \\ \hline
\end{tabular}
\caption{Sieving Table for Case $n=10$, $a=b=c=0$}
\end{table}

For the case $(a,b,c) \neq (0,0,0)$, we conclude, using (\ref{ultNbcnotzero}), that $N>0$ for $\omega(q^{10}-1) \geq 24$. Use of (\ref{finalsieve1bcneqzero}) improves this to $N>0$ for $\omega(q^{10}-1) \geq 13$, as shown in Table 6. The construction of this table, with regards to complementary divisors, is the same as that for Table 2. The possible exceptions here are $q=5$, $7$, and $11$, the last of which satisfies the sieve for $d=2$, with $2^{\omega(e_{i})}-2^{\omega(d)}=2^{\omega(d)}$ for each $i$. The others must be checked directly, and we do so, with an affirmative outcome. Thus $N_{q,10}(a,b,c)>0$.


\begin{table}[htbp]\label{s_table_10_notallzero}
\centering
\begin{tabular}{||c|c|c||}\hline
$\omega(q^{10}-1)$ & $q_{0}$ & $(2/5)$-root of RHS \\ \hline
1 & 1.12 & 3.78 \\ \hline
2 & 1.21 & 5.44 \\ \hline
3 & 1.41 & 7.51 \\ \hline
4 & 1.71 & 9.33 \\ \hline
5 & 2.17 & 11.14 \\ \hline
6 & 2.80 & 12.90 \\ \hline
7 & 3.72 & 14.59 \\ \hline
8 & 5.00 & 16.31 \\ \hline
9 & 6.84 & 18.02 \\ \hline
10 & 9.57 & 19.67 \\ \hline
11 & 13.50 & 21.41 \\ \hline
12 & 19.37 & 23.15 \\ \hline
13 & 28.08 & 24.96 \\ \hline
\end{tabular}
\caption{Sieving Table for Case $n=10$, $(a,b,c) \neq (0,0,0)$}
\end{table}}
\item{$n=9$. First, we consider the case $a=b=c=0$. Observe first that prime $p$ divides $(q^{2}+q+1)$ if and only if $p=3$ or $p \equiv 1 \pmod{6}$, while $p$ divides $(q^{6}+q^{3}+1)$ if and only if $p=3$ or $p \equiv 1 \pmod{18}$ (again, the reader is referred to \cite{koblitz}). Thus, we will only consider prime divisors of the form $p=3$ or $p \equiv 1 \pmod{6}$, as $Q=(q^{2}+q+1)(q^{6}+q^{3}+1)$. With this in hand, we determine from (\ref{ultNzeros}) that $N>0$ for $\omega(Q) \geq 20$. Use of (\ref{finalsieve1allzeros}) improves this to $N>0$ for $\omega(Q) \geq 9$, as Table 7 shows. In Table 7, we use $e_{1}=d=3$ for $\omega=1$, while using $e_{1}=d=3$ and $e_{2}=21$ for $\omega=2$. For $\omega \geq 3$, we use complementary divisors such that $2^{\omega(e_{i})}-2^{\omega(d)}=2^{\omega(d)}$ for each $i$. For $\omega=3$, $d=3$, while for $\omega \geq 4$, we use $d=21$. The possible exceptions are $q=5$ and $q=7$, which are each checked directly to confirm that $N_{5,9}(0,0,0)>0$ and $N_{7,9}(0,0,0)>0$. Thus, $N_{q,9}(0,0,0)>0$.


\begin{table}[htbp]\label{sievetable_9_zeros}
\centering
\begin{tabular}{||c|c|c||}\hline
$\omega(Q)$ & $q_{0}$ & $(2/3)$-root of RHS \\ \hline
1 & 0.68 & 1.90 \\ \hline
2 & 1.20 & 2.61 \\ \hline
3 & 1.83 & 3.13 \\ \hline
4 & 2.75 & 4.16 \\ \hline
5 & 4.33 & 13.87 \\ \hline
6 & 6.89 & 16.56 \\ \hline
7 & 11.11 & 19.17 \\ \hline
8 & 18.67 & 21.65 \\ \hline
9 & 31.67 & 24.08 \\ \hline
\end{tabular}
\caption{Sieving Table for Case $n=9$, $a=b=c=0$}
\end{table}

For the case $(a,b,c) \neq (0,0,0)$, we conclude, using (\ref{ultNbcnotzero}), that $N>0$ for $\omega(q^{9}-1) \geq 30$. Use of (\ref{finalsieve1bcneqzero}) improves this to $N>0$ for $\omega(q^{9}-1) \geq 14$, as shown in Table 8. The construction of this table, with regards to complementary divisors, is the same as that for Table 2. The possible exceptions here are $q=5$, $7$, $11$, $13$, $19$, $23$, and $25$. The first three values are addressed directly via computer (with an affirmative outcome), while the remaining four each satisfy the sieve for $d=2$, with $2^{\omega(e_{i})}-2^{\omega(d)}=2^{\omega(d)}$ for each $i$. Thus $N_{q,9}(a,b,c)>0$.


\begin{table}[htbp]\label{s_table_9_notallzero}
\centering
\begin{tabular}{||c|c|c||}\hline
$\omega(q^{9}-1)$ & $q_{0}$ & square root of RHS \\ \hline
1 & 1.13 & 5.27 \\ \hline
2 & 1.24 & 8.29 \\ \hline
3 & 1.46 & 12.41 \\ \hline
4 & 1.81 & 16.25 \\ \hline
5 & 2.36 & 20.25 \\ \hline
6 & 3.14 & 24.31 \\ \hline
7 & 4.31 & 28.33 \\ \hline
8 & 5.97 & 32.55 \\ \hline
9 & 8.46 & 36.85 \\ \hline
10 & 12.31 & 41.13 \\ \hline
11 & 18.02 & 45.74 \\ \hline
12 & 26.92 & 50.45 \\ \hline
13 & 40.67 & 55.45 \\ \hline
14 & 61.77 & 60.98 \\ \hline
\end{tabular}
\caption{Sieving Table for Case $n=9$, $(a,b,c) \neq (0,0,0)$}
\end{table}}
\end{enumerate}

\begin{Theorem}
$N_{q,n}(a,b,c)>0$ for all prime powers $q=p^{e}$, $p \geq 5$, and for all $n \geq 9$.
\end{Theorem}

\section{The Cases $n=7$, $8$}

We proceed as in the previous section, stating our results as we go.

\begin{enumerate}
\item{$n=8$. From (\ref{ultNzeros}), we determine that $N_{q,8}(0,0,0)>0$ for $\omega(Q) \geq 91$. Use of (\ref{finalsieve1allzeros}) improves this to $N>0$ for $\omega(Q) \geq 18$, as Table 9 shows. The table is built in the same manner as Table 5. The only possible exception are given in Table 10. The values $q=49$, $67$, $73$, $79$, $89$, $109$, $125$, $137$, $173$, $233$ fall to the sieve ($d=2$ for each, with $2^{\omega(e_{i})}-2^{\omega(d)}=2^{\omega(d)}$ for each $i$), while the remaining prime values each succumb to direct verification.


\begin{table}[htbp]\label{sievetable_8_zeros}
\centering
\begin{tabular}{||c|c|c||}\hline
$\omega(Q)$ & $q_{0}$ & RHS \\ \hline
1 & 0.50 & 5.00 \\ \hline
2 & 0.92 & 11.00 \\ \hline
3 & 1.36 & 23.86 \\ \hline
4 & 1.94 & 41.26 \\ \hline
5 & 2.84 & 65.39 \\ \hline
6 & 4.20 & 96.57 \\ \hline
7 & 6.38 & 134.50 \\ \hline
8 & 9.80 & 181.79 \\ \hline
9 & 15.43 & 238.33 \\ \hline
10 & 25.06 & 303.10 \\ \hline
11 & 41.03 & 381.70 \\ \hline
12 & 68.82 & 471.86 \\ \hline
13 & 117.08 & 577.99 \\ \hline
14 & 200.48 & 707.07 \\ \hline
15 & 347.60 & 861.74 \\ \hline
16 & 613.03 & 1043.82 \\ \hline
17 & 1097.75 & 1259.48 \\ \hline
18 & 1975.04 & 1529.57 \\ \hline
\end{tabular}
\caption{Sieving Table for Case $n=8$, $a=b=c=0$}
\end{table}


\begin{table}[htbp]\label{possexp_8_zeros}
\centering
\begin{tabular}{||c|c|c|c|c|c|c|c|c||}\hline
5 & 7 & 11 & 13 & 17 & 19 & 23 & 25 & 29 \\ \hline
31 & 37 & 41 & 43 & 47 & 49 & 53 & 59 & 67 \\ \hline
73 & 79 & 83 & 89 & 109 & 125 & 137 & 173 & 233 \\ \hline
\end{tabular}
\caption{Possible Exceptions for Case $n=8$, $a=b=c=0$}
\end{table}

For the case $(a,b,c) \neq (0,0,0)$, we conclude, using (\ref{ultNbcnotzero}), that $N>0$ for $\omega(q^{8}-1) \geq 45$. Use of (\ref{finalsieve1bcneqzero}) improves this to $N>0$ for $\omega(q^{8}-1) \geq 17$, as shown in Table 11. The construction of this table, with regards to complementary divisors, is the same as that for Table 2. The possible exceptions are listed in Table 12. The values $q=61$, $73$, $79$, $89$ each satisfy the sieve for $d=2$, with $2^{\omega(e_{i})}-2^{\omega(d)}=2^{\omega(d)}$ for each $i$, while the remaining primes are each resolved in the affirmative via computer check.


\begin{table}[htbp]\label{s_table_8_notallzero}
\centering
\begin{tabular}{||c|c|c||}\hline
$\omega(q^{8}-1)$ & $q_{0}$ & $(2/3)$-root of RHS \\ \hline
1 & 1.15 & 9.16 \\ \hline
2 & 1.28 & 16.74 \\ \hline
3 & 1.54 & 28.60 \\ \hline
4 & 1.95 & 40.92 \\ \hline
5 & 2.63 & 54.76 \\ \hline
6 & 3.63 & 69.81 \\ \hline
7 & 5.17 & 85.52 \\ \hline
8 & 7.47 & 102.84 \\ \hline
9 & 11.06 & 121.34 \\ \hline
10 & 16.84 & 140.49 \\ \hline
11 & 25.87 & 161.93 \\ \hline
12 & 40.63 & 184.66 \\ \hline
13 & 64.63 & 209.65 \\ \hline
14 & 103.42 & 238.19 \\ \hline
15 & 167.34 & 270.27 \\ \hline
16 & 274.87 & 305.74 \\ \hline
17 & 457.60 & 345.27 \\ \hline
\end{tabular}
\caption{Sieving Table for Case $n=8$, $(a,b,c) \neq (0,0,0)$}
\end{table}


\begin{table}[htbp]\label{possexp_8_notallzero}
\centering
\begin{tabular}{||c|c|c|c|c|c|c|c||}\hline
5 & 7 & 11 & 13 & 17 & 19 & 23 & 25 \\ \hline 
29 & 31 & 37 & 41 & 43 & 47 & 49 & 53 \\ \hline
59 & 61 & 67 & 73 & 79 & 83 & 89 & \\ \hline
\end{tabular}
\caption{Possible Exceptions for Case $n=8$, $(a,b,c) \neq (0,0,0)$}
\end{table}

}
\item{$n=7$. We consider the all-zero case first. Observe first that prime $p$ divides $Q$ if and only if $p=7$ or $p \equiv 1 \pmod{14}$. Using this, we determine from (\ref{ultNzeros}) that $N>0$ for $\omega(Q) \geq 266$, where we note that the $266$th such prime is $p=13469$. Use of (\ref{finalsieve1allzeros}) dramatically improves this to $N>0$ for $\omega(Q) \geq 10$, as Table 13 shows. In Table 13, we use $e_{1}=d=7$ for $\omega=1$, while for $\omega \geq 2$, we use complementary divisors such that $2^{\omega(e_{i})}-2^{\omega(d)}=2^{\omega(d)}$ for each $i$, and in particular we use $d=1$ for these values of $\omega$, as they produce better results than the method used for other values of $n$. That they produce better results is due to the fact that we are working with primes of a certain form, as opposed to having no restriction on which primes divide $Q$. 

The list of possible exceptions is given in Table 14. The values $q=67$, $125$, $131$, $139$, $223$, $359$, $389$ are eliminated via the sieve, with $d=1$, while the primes less than $q=127$ are eliminated via computer. This leaves us with the possible exceptions $q=25$, $49$, $121$, $169$, $191$, $197$, $199$, $239$, $269$, a total of $9$ possible exceptions.


\begin{table}[htbp]\label{sievetable_7_zeros}
\centering
\begin{tabular}{||c|c|c||}\hline
$\omega(Q)$ & $q_{0}$ & square of RHS \\ \hline
1 & 1.00 & 25.00 \\ \hline
2 & 2.19 & 74.77 \\ \hline
3 & 4.34 & 156.39 \\ \hline
4 & 9.06 & 270.94 \\ \hline
5 & 20.13 & 418.39 \\ \hline
6 & 45.35 & 601.04 \\ \hline
7 & 109.63 & 816.53 \\ \hline
8 & 267.73 & 1067.72 \\ \hline
9 & 667.22 & 1354.61 \\ \hline
10 & 1707.85 & 1676.92 \\ \hline
\end{tabular}
\caption{Sieving Table for Case $n=7$, $a=b=c=0$}
\end{table}


\begin{table}[htbp]\label{possexp_7_zeros}
\centering
\begin{tabular}{||c|c|c|c|c|c|c|c|c|c||}\hline
5 & 7 & 11 & 13 & 17 & 19 & 23 & 25 & 29 & 37 \\ \hline
41 & 43 & 47 & 49 & 53 & 59 & 67 & 71 & 79 & 97 \\ \hline
103 & 107 & 109 & 113 & 121 & 125 & 127 & 131 & 139 & 169 \\ \hline
191 & 197 & 199 & 223 & 239 & 269 & 359 & 389 & & \\ \hline
\end{tabular}
\caption{Possible Exceptions for Case $n=7$, $a=b=c=0$}
\end{table}

For the case $(a,b,c) \neq (0,0,0)$, we use (\ref{ultNbcnotzero}) to conclude that $N>0$ for $\omega(q^{7}-1) \geq 100$. Use of (\ref{finalsieve1bcneqzero}) improves this to $N>0$ for $\omega(q^{7}-1) \geq 21$, as shown in Table 15. The construction of this table, with regards to complementary divisors, is the same as that for Table 2, with the caveat that, for $13 \leq \omega(q^{7}-1) \leq 21$, we use $d=30$ instead of $d=6$, as this gives better results for the RHS values. The $95$ possible exceptions are listed in Table 16. The $42$ values eliminated by the sieve with $d=2$, with $2^{\omega(e_{i})}-2^{\omega(d)}=2^{\omega(d)}$ for each $i$, are given in Table 17. The $38$ primes in Table 16 less than or equal to $179$ were eliminated via computer check, leaving $15$ possible exceptions that range in value from $q=25$ to $q=361$. 


\begin{table}[htbp]\label{s_table_7_notallzero}
\centering
\begin{tabular}{||c|c|c||}\hline
$\omega(q^{7}-1)$ & $q_{0}$ & RHS \\ \hline
1 & 1.17 & 27.63 \\ \hline
2 & 1.32 & 68.15 \\ \hline
3 & 1.63 & 151.74 \\ \hline
4 & 2.15 & 258.90 \\ \hline
5 & 3.02 & 399.61 \\ \hline
6 & 4.36 & 573.90 \\ \hline
7 & 6.54 & 777.08 \\ \hline
8 & 9.96 & 1024.19 \\ \hline
9 & 15.58 & 1312.77 \\ \hline
10 & 25.21 & 1636.69 \\ \hline
11 & 41.17 & 2027.87 \\ \hline
12 & 68.97 & 2473.11 \\ \hline
13 & 117.23 & 2934.03 \\ \hline
14 & 200.63 & 3380.67 \\ \hline
15 & 347.74 & 3866.04 \\ \hline
16 & 613.17 & 4383.96 \\ \hline
17 & 1097.89 & 4935.80 \\ \hline
18 & 1975.19 & 5539.77 \\ \hline
19 & 3601.44 & 6182.56 \\ \hline
20 & 6621.27 & 6874.56 \\ \hline
21 & 12221.66 & 7628.84 \\ \hline
\end{tabular}
\caption{Sieving Table for Case $n=7$, $(a,b,c) \neq (0,0,0)$}
\end{table}


\begin{table}[htbp]\label{possexp_7_notallzero}
\centering
\begin{tabular}{||c|c|c|c|c|c|c|c||}\hline
5 & 7 & 11 & 13 & 17 & 19 & 23 & 25 \\ \hline 
29 & 31 & 37 & 41 & 43 & 47 & 49 & 53 \\ \hline
59 & 61 & 67 & 71 & 73 & 79 & 83 & 89 \\ \hline
97 & 101 & 103 & 107 & 109 & 113 & 121 & 125 \\ \hline
127 & 131 & 137 & 139 & 149 & 151 & 157 & 167 \\ \hline
169 & 173 & 179 & 181 & 191 & 193 & 197 & 199 \\ \hline
211 & 223 & 227 & 229 & 233 & 239 & 241 & 257 \\ \hline
263 & 269 & 271 & 277 & 281 & 283 & 293 & 307 \\ \hline
311 & 313 & 317 & 331 & 337 & 343 & 359 & 361 \\ \hline
367 & 373 & 379 & 389 & 397 & 401 & 409 & 431 \\ \hline
439 & 463 & 491 & 499 & 509 & 529 & 547 & 571 \\ \hline
613 & 625 & 661 & 691 & 727 & 919 & 953 &     \\ \hline
\end{tabular}
\caption{Possible Exceptions for Case $n=7$, $(a,b,c) \neq (0,0,0)$}
\end{table}


\begin{table}[htbp]\label{sieve_7_dis2}
\centering
\begin{tabular}{||c|c|c|c|c|c|c||}\hline
197 & 227 & 229 & 233 & 239 & 263 & 269 \\ \hline
271 & 277 & 281 & 283 & 293 & 307 & 311 \\ \hline
313 & 317 & 337 & 343 & 359 & 367 & 373 \\ \hline
379 & 389 & 397 & 401 & 409 & 431 & 439 \\ \hline
463 & 491 & 499 & 509 & 529 & 547 & 571 \\ \hline
613 & 625 & 661 & 691 & 727 & 919 & 953 \\ \hline
\end{tabular}
\caption{Values for Case $n=7$, $(a,b,c) \neq (0,0,0)$ Eliminated via the Sieve}
\end{table}
}
\end{enumerate}

In closing this section, the author wishes to make an important comment with regards to the direct verification procedures employed for $n=7$, $8$. At the time that he did computational work for these values of $n$, both time and computational resources were (regrettably) limited, more so than in \cite{comi}. Thus he decided that, as resolution of the non-prime $q$ values would consume a great deal more time and resources than could be allowed (see \cite{comi} for a description of how the non-prime $q$ values were handled there), his time could be best spent in eliminating as many prime values as possible from consideration. Given the results of this section, however, it is reasonable to speculate that none of the non-prime $q$ values listed is indeed a genuine exception. Resolution of these values is left to those whose computational resources are sufficient to the task. 



\end{document}